\documentclass{scrartcl}

\usepackage[english]{babel}
\usepackage[utf8]{inputenc}
\usepackage{mathpazo}
\usepackage{microtype}
\usepackage{amsmath}
\usepackage{amsthm}
\usepackage{graphicx}
\usepackage{xcolor}
\usepackage[caption=false]{subfig}
\usepackage[top=1.25in]{geometry}
\usepackage{enumitem}
\usepackage{todonotes}
\usepackage{algorithm}
\usepackage{algpseudocode}
\usepackage[hyphens]{url}
\usepackage[hidelinks]{hyperref}
\usepackage[capitalise,sort&compress]{cleveref}
\usepackage{placeins}

\usepackage{mathrsfs}
\usepackage{comment}
\usepackage{aliascnt}
\usepackage{multirow}
\usepackage{mathtools,amssymb}
\usepackage{booktabs,tabularx,array}
\crefname{appendix}{appendix}{appendices}
\Crefname{appendix}{Appendix}{Appendices}

\providecommand{\ISMFD}{\textsc{ISFD}}     
\providecommand{\LCPG}{\textsc{LCPG}}  
\providecommand{\Clarabel}{\textsc{Clarabel}}
\providecommand{\IPOPT}{\textsc{Ipopt}}

\newtheorem{theorem}{Theorem}[section]

\newaliascnt{lemma}{theorem}
\newtheorem{lemma}[lemma]{Lemma}
\aliascntresetthe{lemma}

\crefname{lemma}{lemma}{lemmas}
\Crefname{lemma}{Lemma}{Lemmas}
\newaliascnt{proposition}{theorem}
\newtheorem{proposition}[proposition]{Proposition}
\aliascntresetthe{proposition}

\crefname{proposition}{proposition}{propositions}
\Crefname{proposition}{Proposition}{Propositions}
\newaliascnt{corollary}{theorem}

\aliascntresetthe{corollary}

\crefname{corollary}{corollary}{corollaries}
\Crefname{corollary}{Corollary}{Corollaries}

\theoremstyle{definition}
\newaliascnt{remark}{theorem}
\newtheorem{remark}[remark]{Remark}
\aliascntresetthe{remark}

\crefname{remark}{remark}{remarks}
\Crefname{remark}{Remark}{Remarks}

\newtheorem{assumption}{Assumption}

\makeatletter
\DeclareOldFontCommand{\sc}{\normalfont\scshape}{\@nomath\sc}
\makeatother

\crefdefaultlabelformat{#2\textup{#1}#3}
\Crefname{assumption}{Assumption}{Assumptions}
\Crefrangeformat{assumption}{Assumptions #3#1#4--#5#2#6}
\crefname{assumption}{assumption}{assumptions}
\Crefname{assumption}{Assumption}{Assumptions}

\crefrangeformat{assumption}
  {assumptions~#3\textup{#1}#4--#5\textup{#2}#6}
\Crefrangeformat{assumption}
  {Assumptions~#3\textup{#1}#4--#5\textup{#2}#6}

\newcommand{\keywords}[1]
{
	{\small
		\textbf{Key words.} {#1}
		\\
	}
}

\renewcommand{\abstract}[1]
{
	{\small
		\textbf{Abstract.} {#1}
		\\
	}
}

\title{A Feasible-Direction-Type Method for Nonconvex Optimization with Infeasible Starts\footnote{\textbf{Funding}: This material is based upon work supported by the National Science Foundation under grant DMS-2410944. The authors also acknowledge support from the Ronald J.\
and Carol T.\ Beerman Early Career Professorship.}}
\author{Danlin Li\thanks{H.\ Milton Stewart School of Industrial and Systems Engineering, Georgia Institute of Technology, Atlanta, Georgia 30332 (\texttt{dli620@gatech.edu}, \texttt{johannes.milz@isye.gatech.edu}).}
\and Johannes Milz\footnotemark[2]
 }
\date{October 28, 2026}
\begin{document}

\maketitle

\abstract{%
We develop an infeasible-start extension of the Topkis--Veinott feasible-direction method for nonconvex composite optimization. At each main iteration, a strongly convex local model is used to obtain objective descent while limiting increases in constraint violation. A separate feasibility-reduction step is invoked when feasibility improvement is needed. Under a feasibility-reduction condition, we establish the standard first-order rate for reducing the Fritz John residual and constraint violation. We then verify that both proximal-gradient and regularized Gauss--Newton methods applied to minimizing the squared constraint-violation function satisfy the required feasibility-reduction property under a Polyak--\L{}ojasiewicz condition. Numerical experiments compare the reduction methods, including an experimental prox-linear variant, and compare the scalar-regularized and curvature-enhanced subproblem formulations.
}


\section{Introduction}
 We consider the nonconvex composite constrained optimization problem:
\begin{align}
\label{eq:infeasible-fsqp}
\min_{\ell\le x\le u}  F(x)\coloneqq f(x)+\psi_0(x)\quad
\mathrm{s.t.}\quad
 c_i(x)\le 0,\quad i\in\mathcal I,
\quad c_i(x)=0,\quad i\in\mathcal E.
\end{align}
Here $x\in\mathbb R^n$, $\mathcal I\coloneqq\{1,\ldots,m\}$, $\mathcal E\coloneqq\{m+1,\ldots,m+p\}$, $\mathcal C\coloneqq\mathcal I\cup\mathcal E$. The bounds $\ell\in(\mathbb R\cup\{-\infty\})^n$ and $u\in(\mathbb R\cup\{+\infty\})^n$ satisfy $\ell\leq u$. The functions $f:\mathbb R^n\to\mathbb R$ and $c_i:\mathbb R^n\to\mathbb R$, $i\in\mathcal C$, are smooth but may be nonconvex. The function $\psi_0:\mathbb R^n\to\mathbb R$ is finite-valued, closed, and convex, allowing the objective to include a structured nonsmooth regularizer. 

Problem \eqref{eq:infeasible-fsqp} combines several features that arise frequently in applications: nonlinear inequality and equality constraints, simple bounds, and a possibly nonsmooth regularizer. The considerations below motivate the design of
the method developed in this work.

A central goal of this work is to develop an algorithm achieving first-order stationarity guarantees without imposing classical constraint qualifications such as the Mangasarian--Fromovitz constraint qualification (MFCQ). Such conditions can fail in important optimization problems, such as mathematical programs with complementarity constraints \cite{FlegelKanzow2003}, and mathematical programs with vanishing constraints \cite{AchtzigerKanzow2008}. In the absence of a constraint qualification, the Karush--Kuhn--Tucker (KKT) conditions need not hold at a local minimizer. We therefore adopt Fritz John optimality conditions, as in the convergence analyses of several feasible-direction methods \cite{TopkisVeinott1967,CawoodKostreva1994,ChenKostreva1999,Birge2000}. 

Constraint degeneracy can also have important numerical consequences. At a local solution with a nonempty KKT multiplier set, the MFCQ is equivalent to boundedness of that set. The homogeneity of the Fritz John conditions allows the multipliers to be normalized, yielding bounded Fritz John multipliers even when constraint qualifications fail. 
The drinking-water network problem in \cite[Section~4.2]{HaeserHinderYe2021} provides a concrete example: its solution satisfies the KKT conditions while MFCQ fails, and severe multiplier growth is reported for \IPOPT{} \cite{WachterBiegler2006}. For this problem, our method recovers approximate KKT multipliers of moderate norm from the Fritz John multipliers, despite the unbounded KKT multiplier set.

A second goal is to allow for infeasible initial and intermediate iterates. Classical feasible-direction methods for smooth inequality-constrained optimization start from a feasible point, compute a feasible descent direction using a linear or quadratic subproblem, and select a step length that preserves feasibility while decreasing the objective \cite{Zoutendijk1960,TopkisVeinott1967,Birge2000}. Our framework permits initial points that satisfy the box constraints but potentially violate functional constraints, even when no equality constraints are present. The main step seeks objective descent without increasing any inequality violation, while a separate reduction step provides feasibility improvement when needed.

Allowing infeasible iterates is particularly relevant in the presence of nonlinear equality constraints. A natural way to account for equality constraints when computing a search direction is to impose a first-order tangency condition. For affine equalities, this condition ensures that the equality residual is preserved exactly. For nonlinear equalities, the residual is preserved only to first order, and nonlinear effects may increase equality violation. This motivates controlling possible increases in nonlinear equality violation during the main step, rather than requiring the violation to decrease at every iteration.

Together, these considerations lead to the Infeasible-Start Feasible-Direction-type method (\ISMFD{}) developed here. The method extends the Topkis--Veinott feasible-direction method \cite{TopkisVeinott1967,Birge2000} to infeasible starts, nonlinear equalities, and composite objectives.

\subsection{Our Contributions}
We introduce \ISMFD{}, an infeasible-start feasible-direction method based on a regularized max model that retains the nonsmooth objective increment, shifts each inequality linearization by its current violation, and imposes a first-order tangency condition for the equalities. The accepted main step combines objective descent with componentwise inequality-violation control and a bound on the increase in equality violation. A separate feasibility-reduction step is invoked when needed.
Our main contributions are as follows.
\begin{enumerate}
\item  
Under a feasibility-reduction condition and  regularity assumptions on the problem data, we establish an \(O(\varepsilon^{-2})\) outer-iteration bound for computing an iterate whose normalized Fritz John proximal residual and constraint violation are simultaneously at most \(\varepsilon\). The result applies to exact solutions of the scalar-regularized main subproblems and does not require MFCQ. The analysis uses a normalized merit function and a uniform upper bound on the regularization parameters selected by the algorithm.

    \item
Under a Polyak--\L{}ojasiewicz (PL) condition for the box-constrained minimization of the squared constraint-violation function, we verify that proximal-gradient (PG) and regularized Gauss--Newton (GN) methods applied to this problem satisfy the feasibility-reduction condition. Both methods uniformly contract the constraint violation. We also show that the abstract feasibility-reduction condition itself implies a feasibility error bound.

\item
We evaluate \ISMFD{} on a subset of the \textsc{CUTEst} test collection \cite{GouldOrbanToint2015CUTEst}, where we also test an experimental prox-linear feasibility-reduction scheme and a curvature-enhanced subproblem variant. We compare \ISMFD{} with the level-constrained proximal-gradient (\LCPG{}) method \cite{Boob2025} on convex and nonconvex quadratically constrained quadratic programs (QCQPs). We also illustrate \ISMFD's multiplier behavior on the drinking-water network problem.
    
\end{enumerate}

\subsection{Related Work}

\paragraph{Feasibility-preserving first-order methods.}
Classical feasible-direction methods for smooth inequality-constrained optimization maintain feasible iterates through suitable search directions and step sizes \cite{Zoutendijk1960,TopkisVeinott1967,Chen2000}. 
For several classical variants, linear and sublinear convergence rates are established in \cite{PironneauPolak1973} under convexity and additional regularity assumptions. Norm-relaxed feasible-direction methods modify the classical direction-finding subproblem by incorporating a matrix-weighted norm \cite{CawoodKostreva1994,ChenKostreva1999}. 
Under weak assumptions, these methods either terminate at a Fritz John point or generate a sequence for which every accumulation point is a Fritz John point. A two-stage feasible-direction method first constructs a descent direction and then modifies it to obtain a feasible one. The paper establishes global convergence of its feasible-direction algorithm and suggests an extension to nonlinear equalities through an inequality reformulation \cite{Herskovits1986}.

For nonsmooth constrained problems, feasibility-preserving methods also use proximal subproblems. The inexact proximal-point method in \cite{Jia2025} 
considers Lipschitz continuous and weakly convex objective and inequality constraint functions, which may be nonsmooth. Starting from a feasible point, the method approximately solves strongly convex proximal subproblems using a switching subgradient method and preserves feasibility of its outer iterates until the target accuracy is reached. It provides near-Fritz John stationarity guarantees without a constraint qualification, whereas its near-KKT guarantee additionally requires a uniform strong MFCQ condition for the proximal subproblems. 
The \LCPG{} method \cite{Boob2025}, which we mentioned above, constructs strongly convex surrogate subproblems through partial linearization and quadratic regularization while gradually increasing constraint levels and maintaining strict feasibility. Its analysis assumes a strictly feasible initial point and uses a uniform MFCQ condition to establish convergence to KKT points. In contrast, \ISMFD{} permits infeasible starts and establishes a normalized Fritz John stationarity guarantee without imposing classical constraint qualifications such as MFCQ.

\paragraph{First-order methods allowing infeasibility.}
In traditional two-phase approaches, Phase~I focuses on finding a feasible point, while Phase~II turns to objective minimization. Combined Phase-I--Phase-II feasible-direction methods integrate these tasks and allow infeasible starts \cite{PolakTrahanMayne1979}. 
Their modified Topkis--Veinott model shifts the objective and constraint branches using the maximum constraint violation, whereas ours shifts only inequality branches by their individual violations.
The unified variant in \cite{PolakHe1991} uses the same direction and step-size rules at feasible and infeasible iterates for semi-infinite optimization. The short-step steepest-descent method in \cite{CartisGouldToint2014} uses separate feasibility and optimization phases, with a corrected worst-case evaluation bound based on scaled first-order stopping conditions \cite{CartisGouldToint2017Corrigendum}.

Other methods control feasibility through the construction and correction of individual steps. Strongly sub-feasible methods allow infeasible starts while preserving inequalities already satisfied, including in problems with nonsmooth objectives \cite{TangJian2012}. The inexact-restoration method in \cite{BuenoMartinez2020} combines contractive restoration with regularized tangential steps and an adaptive merit parameter. Our method instead treats inequalities directly in a composite max model and invokes reduction conditionally
after the main step. Objective-function-free methods adaptively use normal steps to reduce constraint violation and tangential steps to improve optimality, with iteration-complexity guarantees under uniform full-rank Jacobian assumptions \cite{gratton2026,bellavia2026}. The general-constraint extension additionally assumes a PL-type feasibility condition requiring a normal-step criticality measure to uniformly dominate constraint violation \cite[Assumption~AS.8]{bellavia2026}.

For nonsmooth and composite problems, infeasibility is also handled through bundle, regularized, and proximal models. Infeasible proximal bundle methods address nonsmooth convex constrained optimization \cite{Sagastiz2005}. For weakly convex objectives and inequality constraints, quadratically regularized subgradient methods start from an approximately feasible point and solve strongly convex subproblems inexactly using a switching subgradient method \cite{MaLinYang2020}. 
Ghost-penalty frameworks accommodate infeasible starts, nonlinear constraints, and composite objectives \cite{FacchineiEtAl2021Ghost,FacchineiEtAl2022Composite}. Their generalized stationarity notion permits infeasible stationary points. The $O(\varepsilon^{-2})$ complexity bounds are established under extended MFCQ or with a feasible start and convex upper approximations of the constraints \cite{FacchineiEtAl2021Ghost}.
Related proximal-gradient methods use normal--tangential step decompositions \cite{DaiQuRobinson2025,curtis2026}. The KKT-residual complexity bound in \cite{curtis2026} requires a uniform rank condition, whereas our bound controls a normalized Fritz John residual and constraint violation under the feasibility-reduction condition. Augmented-Lagrangian frameworks approximately solve generally nonconvex inner subproblems
\cite{Marchi2023,SahinEtAl2019}, whereas our main step solves a strongly convex max-type model.

\paragraph{Second-order and related methods.}
Feasible sequential quadratic programming methods combine inequality feasibility with objective descent and achieve local two-step superlinear convergence under suitable assumptions \cite{Panier1993}. Related strongly sub-feasible direction methods permit infeasible starts while preserving inequalities already satisfied. Their numerical implementations use quasi-Newton updates to approximate Hessian information \cite{JianZhengHuTang2005,JianZhengTangHu2006}. The strongly sub-feasible framework also handles nonlinear equalities through an auxiliary inequality-constrained formulation \cite{JianXuHan2006}.

Worst-case complexity guarantees have also been established for methods using second- and higher-order models. For equality-constrained problems, a two-phase trust-funnel method uses second-order models and allows objective improvement during its feasibility phase, while retaining worst-case iteration-complexity bounds \cite{CurtisRobinsonSamadi2018}. The evaluation-complexity analysis in \cite{BirginEtAl2016Unscaled} uses higher-order models and measures approximate KKT accuracy without scaling by multiplier norms. A log-barrier method based on trust-region subproblems provides a worst-case iteration bound for obtaining approximate Fritz John points from strictly feasible initial points \cite{HinderYe2024}.

\subsection{Notation}

All vector inequalities and positive parts are understood componentwise. Throughout, $\|\cdot\|_2$ denotes the Euclidean norm for vectors and the spectral norm for matrices. Set \(X\coloneqq[\ell,u]\), and for \(t\in\mathbb R\), let
\(t_+\coloneqq\max\{t,0\}\). For
\(\mathcal A\in\{\mathcal I,\mathcal E,\mathcal C\}\), write
\[
c_{\mathcal A}(x)
\coloneqq
(c_i(x))_{i\in\mathcal A},
\quad \text{and} \quad 
c^+(x)
\coloneqq
\begin{bmatrix}
(c_{\mathcal I}(x))_+\\
c_{\mathcal E}(x)
\end{bmatrix}.
\]

For a differentiable mapping \(G=(G_1,\ldots,G_r):\mathbb R^n\to\mathbb R^r\), we denote its gradient by
$
\nabla G(x) \coloneqq
\begin{bmatrix}
\nabla G_1(x)&\cdots&\nabla G_r(x)
\end{bmatrix}
\in\mathbb R^{n\times r}
$.
For a set \(D \subset \mathbb{R}^n\), \(\delta_D\) denotes its indicator function, and $\operatorname{dist}_2(z,D)\coloneqq\inf_{y\in D}\|z-y\|_2$.
If \(D\) is closed and convex, \(N_D(x)\) denotes its normal cone. For a proper closed convex function \(h\) on $\mathbb{R}^n$, define its prox operator
\[
\operatorname{prox}_h(z)
\coloneqq
\operatorname*{argmin}_{y\in\mathbb R^n}
\big\{
h(y)+(1/2)\|y-z\|_2^2
\big\}.
\]

\section{Algorithmic Framework and Main Result}
\label{sec:framework}
This section develops the \ISMFD{} framework and states its main nonasymptotic guarantee. We state the standing assumptions, define the max-type main-step model and its multipliers, and introduce the normalized Fritz John residuals and the reduction routine. We then prove finite termination of the parameter searches and uniform bounds on the parameters actually selected by the algorithm before stating the complexity theorem.

\subsection{Standing assumptions and derived constants}
\label{subsec:standing-assumptions}

\begin{assumption}[Standing regularity and boundedness]
\label{ass:standing-regularity}
The functions \(f\) and \(c_i\), \(i\in\mathcal C\), are continuously differentiable on an open neighborhood of \(X\). Moreover, the following conditions hold.
\begin{enumerate}[
    label=\textup{(\roman*)},
    ref=\textup{(\roman*)},
    nosep
]

\item
The gradient \(\nabla f\) is \(L_{\nabla f}\)-Lipschitz continuous on
\(X\) for some \(L_{\nabla f}>0\). For every \(i\in\mathcal C\), the gradient \(\nabla c_i\) is \(L_{\nabla c_i}\)-Lipschitz continuous on
\(X\) for some \(L_{\nabla c_i}>0\).

\item
The gradient $\nabla f$ is uniformly bounded on \(X\), that is, 
there exists \(L_f>0\) such that
\(\|\nabla f(x)\|_2\le L_f\) for every \(x\in X\).

\item
The gradients of the inequality and equality constraints are uniformly bounded on \(X\), that is, there exist \(L_{c,\,\mathcal I},L_{c,\,\mathcal E}>0\) such that \(\|\nabla c_i(x)\|_2\le L_{c,\,\mathcal I}\) for every \(x\in X\) and \(i\in\mathcal I\), and \(\|\nabla c_i(x)\|_2\le L_{c,\,\mathcal E}\) for every \(x\in X\) and \(i\in\mathcal E\).

\item
The function \(\psi_0\) is finite-valued, closed, convex, and \(L_{\psi_0}\)-Lipschitz continuous on \(X\) for some \(L_{\psi_0}>0\).
\end{enumerate}
\end{assumption}
We next give sufficient conditions for the standing assumptions.
\begin{remark}
\label{rem:standing-regularity}
Suppose that \(X=[\ell,u]\) is bounded. If \(f\) and \(c_i\), \(i\in\mathcal C\), are twice continuously differentiable on an open neighborhood of \(X\), then compactness of \(X\) gives the uniform bounds in \Cref{ass:standing-regularity}, while boundedness of their Hessians on \(X\) gives the required Lipschitz continuity of their gradients. Moreover, a real-valued convex function \(\psi_0\) is Lipschitz continuous on \(X\). Thus, \Cref{ass:standing-regularity} holds in this setting.
\end{remark}
\begin{assumption}[Boundedness from below]
\label{ass:obj-lower-bound}
There exists \(F_{\inf}\in\mathbb R\) such that
$F(x)\ge F_{\inf}$
for all $x \in X$.
\end{assumption}

\paragraph{Lipschitz constants.}
Under \Cref{ass:standing-regularity}, set
\begin{equation}
\label{eq:aggregate-main-constants}
\bar L_{\mathcal I} \coloneqq \max\Big\{ L_{\nabla f}, \max_{i\in\mathcal I}L_{\nabla c_i} \Big\}, \quad L_{\mathcal E} \coloneqq \frac12\bigg(\sum_{i\in\mathcal E}L_{\nabla c_i}^2\bigg)^{1/2}, \quad
L_{\nabla c} \coloneqq \bigg(\sum_{i\in\mathcal C}L_{\nabla c_i}^2\bigg)^{1/2}.
\end{equation}
The maximum and sum over an empty index set are understood to be zero.

\subsection{Main-step model and ISFD algorithm}
At iteration $k$, let $\tilde x^k\in X$ denote the auxiliary point available from the preceding iteration. We define the max model 
\begin{equation*}
\begin{aligned}
p_k(d)\coloneqq\max\Big\{&\nabla f(\tilde{x}^k)^\top d+\psi_0(\tilde{x}^k+d)-\psi_0(\tilde{x}^k),\\
&c_i(\tilde{x}^k)+\nabla c_i(\tilde{x}^k)^\top d-(c_i(\tilde{x}^k))_+,\quad i\in\mathcal I\Big\}.
\end{aligned}
\end{equation*}
The model combines a first-order approximation of the smooth objective  \(f\) with the exact increment of \(\psi_0\), together with shifted affine approximations of the inequality constraints. 
For a regularization parameter $M_k>0$, we consider the subproblem
\begin{equation}
\label{eq:subproblem'}
\min_{d\in\mathbb R^n}
p_k(d)+(M_k/2)\|d\|_2^2
\quad \text{s.t.} \quad 
\tilde x^k+d\in X,
\quad \nabla c_{\mathcal E}(\tilde x^k)^\top d=0.
\end{equation}
This problem has a unique solution $d^k$. If \(\psi_0\equiv0\), \(\mathcal E=\varnothing\), and \(\tilde x^k\) is feasible, the subproblem \eqref{eq:subproblem'} yields a regularized Topkis--Veinott-type direction model \cite{TopkisVeinott1967,Birge2000}. 

We require the pair $(M_k, d^k)$ to satisfy
\begin{subequations}\label{eq:acceptance}
\begin{align}
F(\tilde{x}^k+d^k)-F(\tilde{x}^k)
&\le -\frac{M_k}{2}\|d^k\|_2^2, \label{eq:acceptance-a}\\
c_{\mathcal{I}}(\tilde{x}^k+d^k) &\le (c_{\mathcal I}(\tilde{x}^k))_+, \label{eq:acceptance-b}\\
\|c_{\mathcal{E}}(\tilde{x}^k+d^k)\|_2 &\le \|c_{\mathcal{E}}(\tilde{x}^k)\|_2+\frac{M_k}{4\tau_k}\|d^k\|_2^2. \label{eq:acceptance-c}
\end{align}
\end{subequations}
\paragraph{Subproblem multipliers.}
We next record the KKT conditions for the subproblem \eqref{eq:subproblem'}. 
We define
\begin{equation}\label{eq: def_h_k}
x^{k+1}\coloneqq \tilde x^k+d^k, \quad \text{and} \quad  h_k\coloneqq p_k(d^k).
\end{equation}
Since $p_k$ is finite-valued and convex and the feasible set of \eqref{eq:subproblem'} is a nonempty polyhedron, the first-order optimality conditions for convex optimization guarantee the existence of multipliers $(\lambda_0^k,\lambda_{\mathcal I}^k,\lambda_{\mathcal E}^k) \in\mathbb R^{m+p+1}$ satisfying
\begin{equation}\label{eq:inf-optconds}
\begin{aligned}
s_0^k \in \partial \psi_0(x^{k+1}), \quad
s_{-1}^k \in N_{X}(x^{k+1}), \quad x^{k+1} \in X,\quad \lambda_j^k&\in \mathbb{R}, \quad j\in \mathcal{E}, 
 \\\lambda^k_i \ge 0, \quad i\in \{0\}\cup\mathcal{I}, \quad
\lambda_0^k+\sum_{i\in\mathcal{I}} \lambda_i^k = 1,\quad
\nabla c_j(\tilde x^k)^{\top}d^k &=0, \quad j\in \mathcal{E}, \qquad
 \\
\lambda_0^k \Bigl(
\nabla f(\tilde{x}^k)^\top d^k
+ \psi_0(x^{k+1})
- \psi_0(\tilde{x}^k)
- h_k
\Bigr)
&= 0, \\
\lambda_{i}^k \Bigl(
c_i(\tilde{x}^k)
+ \nabla c_i(\tilde{x}^k)^\top d^k
- (c_i(\tilde{x}^k))_+
- h_k
\Bigr)
&= 0, \quad i \in \mathcal I, \\
\lambda_0^k\bigl[\nabla f(\tilde{x}^k)+s_0^k\bigr]
+ \sum_{i\in\mathcal{I}} \lambda_i^k \nabla c_i(\tilde{x}^k) +\sum_{j\in \mathcal E}\lambda_j^k\nabla c_j(\tilde x^k)
+ s_{-1}^k
&= -M_k d^k.
\end{aligned}
\end{equation}
In particular, if $d^k=0$, $\lambda_0^k >0$, and $\tilde{x}^k$ is feasible, then $x^{k+1}=\tilde{x}^k$ satisfies the KKT conditions of \eqref{eq:infeasible-fsqp} with multipliers $\lambda_i^k/\lambda_0^k$, $i\in\mathcal C$.
\paragraph{Multiplier normalization and residuals.}
The multipliers in \eqref{eq:inf-optconds} associated with the objective and inequality constraints satisfy   $\lambda_0^k+\sum_{i\in\mathcal I}\lambda_i^k=1$, whereas the equality constraint multipliers are not normalized. Let $\lambda^k \coloneqq (\lambda_{\mathcal I}^k,\lambda_{\mathcal E}^k)$ and define the normalized multipliers
\begin{equation}
\label{eq:multiplier-normalization}
(\bar\lambda_0^k,\bar\lambda^k) \coloneqq (\lambda_0^k,\lambda^k)/\kappa_k,
\quad \text{where} \quad \kappa_k \coloneqq 1+\|\lambda_{\mathcal E}^k\|_1.
\end{equation}
It follows that
\begin{equation}
    \bar\lambda_0^k +\|\bar\lambda_{\mathcal{I}}^k\|_1 +\|\bar\lambda_{\mathcal{E}}^k\|_1 = 1.
\end{equation}
Since $\kappa_k>0$, the normalization preserves the sign and complementarity relations in \eqref{eq:inf-optconds}.

Motivated by the optimality system \eqref{eq:inf-optconds} and the Fritz John conditions for \eqref{eq:infeasible-fsqp}, for $x\in X$, $\lambda_0\ge 0$, and $\lambda\in\mathbb{R}^{m+p}$ with $\lambda_{\mathcal I}\ge 0$, define $g(x,\lambda_0,\lambda)$ by
\[
g(x,\lambda_0,\lambda) \coloneqq \lambda_0\nabla f(x) + \sum_{i\in\mathcal C}\lambda_i\nabla c_i(x),
\]
and the following two Fritz John residuals:
\begin{align}
r_{\mathrm{prox}}(x,\lambda_0,\lambda)
\coloneqq
\max\Big\{&
\big\|
x-
\operatorname{prox}_{\lambda_0\psi_0+\delta_X}
\bigl(x-g(x,\lambda_0,\lambda)\bigr)
\big\|_2,
\;
\sum_{i\in\mathcal I}\lambda_i|c_i(x)|
\Big\},
\label{eq:residual-r-prox}
\\
r_{\mathrm{dist}}(x,\lambda_0,\lambda)
\coloneqq
\max\Big\{&
\operatorname{dist}_2\!\left(
-g(x,\lambda_0,\lambda),
\lambda_0\partial\psi_0(x)+N_X(x)
\right),
\;
\sum_{i\in\mathcal I}\lambda_i|c_i(x)|
\Big\}.
\label{eq:residual-r-dist}
\end{align}
For $\lambda_{\mathcal{I}} \geq 0$, 
their common second term measures complementarity, while infeasibility is measured separately by $\|c^+(x)\|_2$. At iteration $k$, the termination test and reduction trigger of \Cref{alg:infeasiblefsqp} use $r_{\mathrm{prox}}(x^{k+1},\bar\lambda_0^k,\bar\lambda^k)$, whereas the convergence analysis uses $r_{\mathrm{dist}}(x^{k+1},\bar\lambda_0^k,\bar\lambda^k)$. The following lemma shows that the Fritz John distance residual dominates the Fritz John proximal residual and that one vanishes if and only if the other does.
\begin{lemma}[Residual comparison]
\label{lem:prox-distance-relation}
Suppose that \Cref{ass:standing-regularity} holds. For every \(x\in X\), \(\lambda_0\ge0\), and
\(\lambda=(\lambda_{\mathcal I},\lambda_{\mathcal E})\) with
\(\lambda_{\mathcal I}\ge0\),
we have
\[
r_{\mathrm{prox}}(x,\lambda_0,\lambda)
\le
r_{\mathrm{dist}}(x,\lambda_0,\lambda).
\]
Moreover,
$r_{\mathrm{prox}}(x,\lambda_0,\lambda)=0$
if and only if 
$r_{\mathrm{dist}}(x,\lambda_0,\lambda)=0$.
\end{lemma}

\begin{proof}
Using calculus rules for subdifferentials and applying \cite[Lemma~4.1.6]{Milzarek2016NumericalMA}  yield
\[
\big\|x-\operatorname{prox}_{\lambda_0\psi_0+\delta_X}\bigl(x-g(x,\lambda_0,\lambda)\bigr)\big\|_2\le \operatorname{dist}_2(-g(x,\lambda_0,\lambda), \lambda_0\partial\psi_0(x)+N_X(x)). 
\]
Hence, we obtain the inequality. If \(r_{\mathrm{prox}}(x,\lambda_0,\lambda)=0\), then optimality conditions ensure $-g(x,\lambda_0,\lambda)\in\lambda_0\partial\psi_0(x)+N_X(x)$, and, hence, \(r_{\mathrm{dist}}(x,\lambda_0,\lambda)=0\). 
\end{proof}
Consequently, for any normalized multiplier pair \((\bar\lambda_0,\bar\lambda)\) satisfying \eqref{eq:multiplier-normalization}, \(\bar\lambda_0\ge0\) and \(\bar\lambda_{\mathcal I}\ge0\),
the conditions
\[
r_{\mathrm{prox}}(x,\bar\lambda_0,\bar\lambda)=0, \quad  \text{and} \quad  \|c^+(x)\|_2=0,
\]
are equivalent to the Fritz John stationarity, complementarity, and feasibility conditions.

\paragraph{Reduction routine.}
We introduce a map for reducing constraint violation. The abstract
analysis requires only a displacement bound and a uniform contraction
of the constraint violation. In \Cref{sec:reduction-step}, we show that
the proximal-gradient (PG) and regularized Gauss--Newton (GN) schemes
satisfy these properties under suitable regularity conditions.

\begin{assumption}[Contractive reduction map]
\label{ass:reduction-contraction}
Let $\theta>0$. There exist a map
$\mathcal R_\theta:X\to X$ and a constant
$\eta_{\mathcal R_\theta}\in(0,1)$
such
that, for every $x\in X$, the point
$\widehat x\coloneqq\mathcal R_\theta(x)\in X$ satisfies
\[
\|\widehat x-x\|_2 \le \theta\|c^+(x)\|_2, 
\quad \text{and} \quad 
\|c^+(\widehat x)\|_2 \le \eta_{\mathcal R_\theta}\|c^+(x)\|_2.
\]
\end{assumption}

Let
\(
S\coloneqq\{x\in X:c_{\mathcal I}(x)\le0,\ c_{\mathcal E}(x)=0\}
\)
denote the feasible set of \eqref{eq:infeasible-fsqp}.

\begin{proposition}[Feasibility error bound implied by reduction]
\label{prop:reduction-error-bound}
If \Cref{ass:reduction-contraction} holds, then $S$ is nonempty, and
\[
\operatorname{dist}_2(x,S)\le \theta(1-\eta_{\mathcal{R}_\theta})^{-1}\|c^+(x)\|_2
\quad \text{for all} \quad  x\in X.
\]
\end{proposition}
\begin{proof}
Fix $z^0=x\in X$ and define $z^{j+1}=\mathcal{R}_\theta(z^j)$ for $j\ge0$. \Cref{ass:reduction-contraction}  gives
\[
\|c^+(z^j)\|_2\le\eta_{\mathcal{R}_\theta}^j\|c^+(x)\|_2,
\quad \text{and} \quad 
\|z^{j+1}-z^j\|_2\le\theta\eta_{\mathcal{R}_\theta}^j\|c^+(x)\|_2.
\]
Thus $(z^j)$ is Cauchy and converges to $\bar{z}\in X$. Continuity of $c^+$ implies $\bar{z}\in S$. Summing the displacement bounds yields
$\|\bar{z}-x\|_2\le\theta\|c^+(x)\|_2/(1-\eta_{\mathcal{R}_\theta})$.
\end{proof}
\paragraph{Algorithmic interpretation.}
We are now ready to state the complete \ISMFD{} framework. At iteration \(k\), the method first computes an accepted main point $x^{k+1}=\tilde x^k+d^k$
from the current model center $\tilde x^k$. The normalized Fritz John proximal residual measures stationarity and inequality complementarity, while \(\|c^+(x^{k+1})\|_2\) measures feasibility. If both quantities vanish, the algorithm terminates. If not, the triggering rule determines whether the main point is sufficiently feasible relative to its normalized Fritz John proximal residual. When the triggering condition is satisfied, \(x^{k+1}\) becomes the next model center. Otherwise, the reduction routine gives $ \tilde x^{k+1}=\mathcal{R}_\theta(x^{k+1})$ and increases the merit parameter until \eqref{eq:reduction-step-merit} holds.

\begin{algorithm}[!htbp]
\small
\caption{\ISMFD{} with adaptive merit parameter}
\label{alg:infeasiblefsqp}
    \begin{enumerate}
        \item[0.] Choose $x^{0} \in X$, 
        and constants $\tau_0>0$ and $\underline{M}>0$.
Let $\beta\in (0, 1)$, $\gamma_\tau>1$, and $\theta>0$. Set $\tilde{x}^0 = x^0$, and choose a reduction map $\mathcal{R}_\theta$.
Let $\overline M$ be given by \eqref{eq:overline-M}.
    \end{enumerate}
    For $k=0, 1, \ldots$:
    \begin{enumerate}
    
        \item \label{step:main} Compute $M_k \in [\underline M, \overline M]$ 
        and a solution $d^k$ to \eqref{eq:subproblem'} 
        such that \eqref{eq:acceptance-a}--\eqref{eq:acceptance-c}
        are satisfied, using a first-acceptable doubling search initialized at $\underline M$. Normalize the associated multipliers satisfying \eqref{eq:inf-optconds} according to \eqref{eq:multiplier-normalization}.

\item \label{step:termination} Set $x^{k+1} = \tilde{x}^k + d^k$. If
$r_{\mathrm{prox}}(x^{k+1}, \bar\lambda_0^k, \bar\lambda^k) = 0$
and $\|c^+(x^{k+1})\|_{2} = 0$, then terminate and return $x^{k+1}$ and $(\bar\lambda_0^k, \bar\lambda^k)$.

        \item \label{step:reduction} If
        \begin{equation}   
       \label{eq:triggering-rule}
             \|c^+(x^{k+1})\|_2
            \leq  \beta 
            r_{\mathrm{prox}}(x^{k+1}, \bar\lambda_0^k, \bar\lambda^k),
         \end{equation}
        set
\[
\tilde x^{k+1}=x^{k+1},
\quad \text{and} \quad 
\tau_{k+1}=\tau_k.
\]
        Otherwise, apply the reduction routine $\mathcal{R}_\theta$ and set
        \[
        \tilde x^{k+1} = \mathcal{R}_\theta(x^{k+1}).
        \]
        Then choose the smallest $\tau_{k+1}\in \{\tau_k, \gamma_\tau\tau_k , \gamma_\tau^2\tau_k,\ldots\}$ such that
        \begin{equation}\label{eq:reduction-step-merit}
            F(\tilde{x}^{k+1}) +\tau_{k+1}\|c^+(\tilde{x}^{k+1})\|_2 \le F(x^{k+1})+\tau_{k+1}\|c^+(x^{k+1})\|_2.
            \end{equation}
    \end{enumerate}
\end{algorithm}

\subsection{Main-step well-definedness and uniform parameter bounds}\label{sec:subproblem-property} 
We derive basic estimates for the main-step model and show that the acceptance conditions hold for sufficiently large regularization parameters. We then establish finite termination of the merit-parameter update and prove that the first-acceptable doubling search selects $M_k\in[\underline M,\overline M]$ at every iteration.

Recall that $d^k$ denotes the unique solution of \eqref{eq:subproblem'} and $h_k$ is defined in \eqref{eq: def_h_k}.

\begin{lemma}[Basic main-step estimates]
\label{lem:basicproperties}
Suppose that \Cref{ass:standing-regularity} holds. Fix $k\geq 0$ and $M_k>0$. 
Then the following statements hold.
\begin{enumerate}[label=\textup{(\roman*)}, ref=\textup{(\roman*)}, nosep]
\item
\label{itm:model-decrease}
The max-model value satisfies
\begin{equation}
\label{eq:decrease}
h_k\le-M_k\|d^k\|_2^2.
\end{equation}
\item
\label{itm:descent}
The objective and inequality constraints satisfy
\begin{align*}
F(\tilde x^k+d^k) - F(\tilde{x}^k) &\leq -\bigl( M_k-(L_{\nabla f}/2)\bigr) \|d^k\|_2^2, \quad \text{and}
\\
c_i(\tilde x^k+d^k) - (c_i(\tilde{x}^k))_+ &\leq - \bigl(M_k -(L_{\nabla c_i}/2)\bigr)\|d^k\|_2^2,
\quad   i \in \mathcal{I}.
\end{align*}
\item \label{itm:equality-bound} The equality constraint violation satisfies
\begin{equation}
\label{eq:basic-eq-norm-bound}
\|c_{\mathcal E}(\tilde x^k+d^k)\|_2
\le
\|c_{\mathcal E}(\tilde x^k)\|_2
+
L_{\mathcal E}\|d^k\|_2^2.
\end{equation}
\end{enumerate}
\end{lemma}
\begin{proof}
We define \[
\mathcal{D}_k\coloneqq\{d\in\mathbb{R}^n:\ \ell\le \tilde{x}^k+d\le u, \;\nabla c_{\mathcal{E}}(\tilde x^k)^{\top}d=0\}
\]
and $\widehat p_k(d)\coloneqq p_k(d)+\delta_{\mathcal{D}_k}(d)$.
Then the subproblem \eqref{eq:subproblem'} is equivalent to
\[
\min_{d\in\mathbb{R}^n}\ \widehat p_k(d)+(M_k/2)\|d\|_2^2.
\]
Since $p_k$ is proper, closed and convex, and $\mathcal{D}_k$ is a nonempty closed convex set, $\widehat p_k$ is also proper, closed and convex. Moreover, $0\in \mathcal{D}_k$ because
$\tilde{x}^k\in X$ and $\nabla c_{\mathcal{E}}(\tilde x^k)^{\top}0 =0$. Applying \cite[Lemma 1]{Boob2025}, we obtain
\[
\widehat p_k(d^k)\le \widehat p_k(0)-M_k\|d^k\|_2^2.
\]
Because \(0,d^k\in \mathcal{D}_k\) and \(p_k(0)=0\), this proves part~\ref{itm:model-decrease}.

Since $\nabla f$ is $L_{\nabla f}$-Lipschitz continuous, the descent lemma gives
\begin{equation}
    \label{eq:f-ratio-bound}
\begin{aligned}
    F(\tilde{x}^k+d^k)&\leq F(\tilde{x}^k)+p_k(d^k)+(L_{\nabla f}/2)\|d^k\|^2_2\\
&\overset{\eqref{eq:decrease}}\leq F(\tilde{x}^k)-(M_k-(L_{\nabla f}/2))\|d^k\|^2_2.
\end{aligned}
\end{equation}
Similarly, since each $\nabla c_i$ is $L_{\nabla c_i}$-Lipschitz continuous and \eqref{eq:decrease}, we have
\begin{align*}
    c_i(\tilde{x}^k+d^k) - (c_i(\tilde{x}^k))_+ \leq -(M_k-(L_{\nabla c_i}/2))\|d^k\|^2_2,
\end{align*}
for all $ i\in \mathcal{I}$.
This proves part \ref{itm:descent}.

For part~\ref{itm:equality-bound}, the componentwise Lipschitz-gradient Taylor estimate and $\nabla c_{\mathcal E}(\tilde x^k)^\top d^k=0$ give
\[
\|c_{\mathcal E}(\tilde x^k+d^k) -c_{\mathcal E}(\tilde x^k)\|_2\le L_{\mathcal E}\|d^k\|_2^2.
\]
The triangle inequality then yields \eqref{eq:basic-eq-norm-bound}.
\end{proof}
For the parameters of \Cref{alg:infeasiblefsqp}, set \begin{equation}
\label{eq:tau-choice}
\tau_{\Theta}
\coloneqq \frac{\theta (L_f+L_{\psi_0})} {1-\eta_{\mathcal{R}_\theta}}, \quad \text{and} \quad  \bar\tau \coloneqq \max\{\tau_0,\gamma_\tau\tau_{\Theta}\}.
\end{equation}
and let
\begin{align}
\label{eq:overline-M}
\overline M\coloneqq 2\max\{\underline M,M(\bar\tau)\},
\quad \text{where} \quad 
M(\tau)
\coloneqq
\max\big\{
L_{\nabla f},\,
(1/2)\max_{i\in\mathcal I}L_{\nabla c_i},\,
4\tau L_{\mathcal E}
\big\}.
\end{align}

\begin{lemma}[Finite termination and uniform parameter bounds]\label{lem:parameter-searches}
Suppose that \Cref{ass:standing-regularity,ass:reduction-contraction} hold. Let $\tau_0>0$ and $\gamma_\tau>1$. Then, for every performed iteration $k\ge0$ of \Cref{alg:infeasiblefsqp}, the following statements hold.  
\begin{enumerate}[label=\textup{(\roman*)}, ref=\textup{(\roman*)}, nosep]
\item \label{itm:tau-search} If the reduction routine is invoked at iteration $k$, then the subsequent merit-parameter search in Step~\ref{step:reduction} terminates finitely, and the resulting merit parameter satisfies $\tau_{k+1}\le\bar\tau$. If Step~\ref{step:reduction} is reached and the reduction routine is not invoked, then $\tau_{k+1}=\tau_k\le\bar\tau$.
\item \label{itm:sufficient-M} If \(M_k\ge M(\tau_k)\),
then $d^k$ satisfies 
\eqref{eq:acceptance-a}--\eqref{eq:acceptance-c}.
\item \label{itm:M-search} The first-acceptable doubling search in Step~\ref{step:main} terminates finitely, and its selected parameter satisfies
\(\underline M\le M_k\le\overline M \).
\end{enumerate}
\end{lemma}
\begin{proof}
For part~\ref{itm:tau-search}, suppose first that the reduction routine is invoked at iteration \(k\). Since \(F\) is \((L_f+L_{\psi_0})\)-Lipschitz on \(X\), \Cref{ass:reduction-contraction} gives, for every \(\tau>0\),
\[
\begin{aligned}
&F(\tilde x^{k+1})+\tau\|c^+(\tilde x^{k+1})\|_2 -
\big[F(x^{k+1})+\tau\|c^+(x^{k+1})\|_2 \big] \\
&\qquad\le \bigl[\theta (L_f+L_{\psi_0})-\tau(1-\eta_{\mathcal R_{\theta}})\bigr]\|c^+(x^{k+1})\|_2.
\end{aligned}
\]
Hence, every \(\tau\ge\tau_\Theta\) satisfies \eqref{eq:reduction-step-merit}, so the geometric search terminates finitely and
\[
\tau_{k+1}\le \max\{\tau_k,\gamma_\tau\tau_\Theta\}.
\]
If the reduction routine is not invoked, then \(\tau_{k+1}=\tau_k\), so the same bound holds. Since \(\bar\tau=\max\{\tau_0,\gamma_\tau\tau_\Theta\}\), induction yields \(\tau_k\le\bar\tau\) for all \(k\ge0\).

For part~\ref{itm:sufficient-M}, \Cref{lem:basicproperties}~\ref{itm:descent} and the definition of \(M(\tau_k)\) give \eqref{eq:acceptance-a}--\eqref{eq:acceptance-b}. Similarly, \Cref{lem:basicproperties}~\ref{itm:equality-bound} and \(M_k\ge4\tau_kL_{\mathcal E}\) yield \eqref{eq:acceptance-c}.

For part~\ref{itm:M-search}, part~\ref{itm:tau-search} gives \(\tau_k\le\bar\tau\), while part~\ref{itm:sufficient-M} shows that every trial \(M\ge M(\tau_k)\) is acceptable. Hence the doubling search terminates finitely. If the initial value \(\underline M\) is accepted, then \(M_k=\underline M\). Otherwise, the trial immediately preceding \(M_k\) was rejected, so part~\ref{itm:sufficient-M} implies \(M_k/2<M(\tau_k)\). Therefore,
\[
\underline M \le M_k \le \max\{\underline M,2M(\tau_k)\} \le 2\max\{\underline M, M(\bar\tau)\}=\overline M,
\]
where the last inequality follows from the monotonicity of \(M(\cdot)\).
\end{proof}

\subsection{Main nonasymptotic result}\label{sec:main-complex-result}
By \Cref{lem:parameter-searches}, the merit-parameter backtracking and the main-step regularization search terminate finitely, and every performed iteration satisfies \(\tau_k\le\bar\tau\) and \(\underline M\le M_k\le\overline M\). We now collect the constants used in the complexity bound and state the main result.
Define
\begin{equation}\label{eq:Cmu-Cd-def}
C_1 \coloneqq \overline{M}+\bar L_{\mathcal I}+2L_{\mathcal E}, \quad \text{and} \quad 
C_2 \coloneqq \bar L_{\mathcal I}+L_{\mathcal{E}}.
\end{equation}
Next, define
\begin{equation}\label{eq:Delta-D-def}
\Delta_{\Theta}\coloneqq \tau_0^{-1}(F(x^0)-F_{\inf})+\|c^+(x^0)\|_2, \quad \text{and} \quad  D_{\Theta}\coloneqq\frac{4\bar\tau\Delta_{\Theta}}{\underline M}.
\end{equation}
Let
\[
Q_\Theta
\coloneqq
\big(
C_1+
\left(1+C_2/\underline M\right)
(L_f+L_{\psi_0})
\big)^2D_\Theta.
\]
Finally, set
\[
\rho\coloneqq \max\{\eta_{\mathcal{R}_\theta},\,\beta\}.
\]
\begin{theorem}[Iteration complexity]
\label{thm:inf-nonasymp}
Suppose that \Cref{ass:standing-regularity,ass:obj-lower-bound,ass:reduction-contraction} hold. Let $K\ge 0$. If \Cref{alg:infeasiblefsqp} does not terminate during iterations $0,\ldots, K$, then

\begin{equation}
\label{eq:main-nonasymptotic-bound}
\min_{0\leq k\leq K}
\max\Big\{
r_{\mathrm{dist}}(x^{k+1},\bar\lambda_0^k,\bar\lambda^k),
\|c^+(x^{k+1})\|_2
\Big\}
\leq
\frac{1}{1-\rho}
\left(
\frac{\|c^+(x^0)\|_2}{K+1}
+\sqrt{\frac{Q_\Theta}{K+1}}
\right).
\end{equation}

The same bound holds with \(r_{\mathrm{prox}}(x^{k+1},\bar\lambda_0^k,\bar\lambda^k)\) in place of \(r_{\mathrm{dist}}(x^{k+1},\bar\lambda_0^k,\bar\lambda^k)\). 
\end{theorem}
The proof of \Cref{thm:inf-nonasymp} is deferred to \Cref{sec:proof-main-theorem}. 
Since the right-hand side of \eqref{eq:main-nonasymptotic-bound} is \(O((K+1)^{-1/2})\), for every \(\varepsilon>0\), within \(O(\varepsilon^{-2})\) outer iterations the algorithm either terminates at an exact feasible Fritz John point or generates an iterate \(x^{k+1}\) satisfying
\[
r_{\mathrm{prox}} \bigl(x^{k+1},\bar\lambda_0^k,\bar\lambda^k\bigr) \le \varepsilon, \quad \text{and} \quad \|c^+(x^{k+1})\|_2\le \varepsilon.
\]
The distance-residual estimate is stronger and is stated as a best-iterate bound over a nonterminated prefix. Although \(r_{\mathrm{prox}}\le r_{\mathrm{dist}}\), small proximal residual alone does not imply a comparably small distance residual. A corresponding KKT complexity statement would require, for example, a uniform positive lower bound on the normalized objective multiplier, which is not assumed here. In \Cref{sec:reduction-step}, we verify the reduction property for PG and GN under additional constraint-value bounds and a uniform feasibility--PL condition.

\section{Convergence and Complexity Analysis}

In this section, we prove \Cref{thm:inf-nonasymp}. A normalized merit function yields square summability of the main steps. The main-step length bound controls the objective-decrease term in the residual estimate. A common recursion for feasibility and stationarity then gives the nonasymptotic result.

\subsection{Merit decrease and summability}
For $\tau>0$, define the merit function
\[
\Theta_\tau(x)\coloneqq
\tau^{-1}(F(x)-F_{\inf})+\|c^+(x)\|_2+\delta_X(x).
\]
\begin{lemma}[Merit decrease and square summability]\label{lem:merit-decrease}
Suppose that \Cref{ass:standing-regularity,ass:obj-lower-bound,ass:reduction-contraction} hold and that \Cref{alg:infeasiblefsqp} does not terminate during iterations $0,\ldots,K$. Then, for $0\le k\le K$,
\begin{equation}\label{eq:merit-decrease}
\begin{aligned}
\Theta_{\tau_{k+1}}(\tilde x^{k+1}) \le\Theta_{\tau_k}(\tilde x^k)-\frac{M_k}{4\tau_k}\|d^k\|_2^2\le\Theta_{\tau_k}(\tilde x^k)-\frac{\underline M}{4\bar\tau}\|d^k\|_2^2.
\end{aligned}
\end{equation}
Consequently,
\begin{equation}\label{eq:sum-dk}
\sum_{k=0}^K\|d^k\|_2^2\le D_\Theta.
\end{equation}
\end{lemma}
\begin{proof}
The inequality and equality acceptance tests, together with
$\sqrt{s^2+(t+\delta)^2}\le\sqrt{s^2+t^2}+\delta$ for $s,t,\delta\ge0$, give
\begin{equation}\label{eq:bounded-increasing-c}
\|c^+(x^{k+1})\|_2
\le\|c^+(\tilde x^k)\|_2+\frac{M_k}{4\tau_k}\|d^k\|_2^2.
\end{equation}
Combining this estimate with \eqref{eq:acceptance-a} and dividing the objective decrease by $\tau_k$ yields
\begin{equation}\label{eq:step1-merit}
\Theta_{\tau_k}(x^{k+1})
\le\Theta_{\tau_k}(\tilde x^k)-\frac{M_k}{4\tau_k}\|d^k\|_2^2.
\end{equation}
Since $\tau_{k+1}\ge\tau_k$ and $F(x^{k+1})\ge F_{\inf}$,
\[
\Theta_{\tau_{k+1}}(x^{k+1})\le\Theta_{\tau_k}(x^{k+1}).
\]
If reduction is skipped, $\tilde x^{k+1}=x^{k+1}$. Otherwise, dividing \eqref{eq:reduction-step-merit} by $\tau_{k+1}$ and subtracting $F_{\inf}/\tau_{k+1}$ gives
\[
\Theta_{\tau_{k+1}}(\tilde x^{k+1})\le\Theta_{\tau_{k+1}}(x^{k+1}).
\]
Thus \eqref{eq:merit-decrease} follows, using $M_k\ge\underline M$ and $\tau_k\le\bar\tau$. Summing it and using the definition of \(\Delta_\Theta\) and the nonnegativity of \(\Theta\), we obtain
\[
\frac{\underline M}{4\bar\tau}\sum_{k=0}^K\|d^k\|_2^2
\le\Theta_{\tau_0}(x^0)-\Theta_{\tau_{K+1}}(\tilde x^{K+1})
\le\Delta_\Theta,
\]
which proves \eqref{eq:sum-dk}.
\end{proof}

\subsection{Residual and feasibility estimates}
\label{subsec:residual-feasibility-estimates}
We next establish the key estimates used in the main complexity analysis. Recall the Fritz John proximal and distance residuals \(r_{\mathrm{prox}}\) and \(r_{\mathrm{dist}}\) defined in \eqref{eq:residual-r-prox} and \eqref{eq:residual-r-dist}. For iterates generated by \Cref{alg:infeasiblefsqp}, define
\[
r^{k+1}_{\mathrm{prox}} \coloneqq r_{\mathrm{prox}}(x^{k+1},\bar\lambda_0^k,\bar\lambda^k),
\quad \text{and} \quad 
r^{k+1}_{\mathrm{dist}} \coloneqq r_{\mathrm{dist}}(x^{k+1},\bar\lambda_0^k,\bar\lambda^k).
\]

\begin{lemma}[Normalized residual estimate]\label{lem:R+bounded}
Suppose that \Cref{ass:standing-regularity} holds. Let $(M_k,d^k)$ be an accepted main step of \Cref{alg:infeasiblefsqp} with $M_k\le\overline M$, let $x^{k+1}=\tilde x^k+d^k$, and let $(\bar\lambda_0^k,\bar\lambda^k)$ be the associated normalized multipliers defined by
\eqref{eq:multiplier-normalization}. Then
\begin{equation}\label{eq:Rk+bound}
\begin{aligned}
r^{k+1}_{\mathrm{prox}}\le r^{k+1}_{\mathrm{dist}}
\le C_1\|d^k\|_2 + C_2\|d^k\|_2^2 + F(\tilde{x}^k) - F(x^{k+1}) +\|c^+(\tilde{x}^k)\|_2,
\end{aligned}
\end{equation}
where $C_1$ and $C_2$ are defined in
\eqref{eq:Cmu-Cd-def}.

\end{lemma}
\begin{proof}
We first bound the stationarity term. Dividing the stationarity condition in \eqref{eq:inf-optconds} by
$\kappa_k=1+\|\lambda_{\mathcal E}^k\|_1$ gives
\[
\bar\lambda_0^k
\bigl( \nabla f(\tilde x^k)+s_0^k\bigr)+\sum_{i\in\mathcal C}
\bar\lambda_i^k\nabla c_i(\tilde x^k)+
\frac{s_{-1}^k}{\kappa_k}=-\frac{M_k}{\kappa_k}d^k.
\]
Since $N_X(x^{k+1})$ is a cone, \[
\bar\lambda_0^ks_0^k + \frac{s_{-1}^k}{\kappa_k} \in \bar\lambda_0^k\partial\psi_0(x^{k+1})+N_X(x^{k+1}).
\]
Using this element in the distance definition and
\(x^{k+1}-\tilde x^k=d^k\), we obtain
\[
\begin{aligned}
&
\operatorname{dist}_2\!\big(
-\bar\lambda_0^k\nabla f(x^{k+1})
-\sum_{i\in\mathcal C}
\bar\lambda_i^k\nabla c_i(x^{k+1}),
\bar\lambda_0^k\partial\psi_0(x^{k+1})
+
N_X(x^{k+1})
\big)
\\
&\qquad\le
\bigg(
\frac{M_k}{\kappa_k}
+
\bar\lambda_0^kL_{\nabla f}
+
\sum_{i\in\mathcal I}
\bar\lambda_i^kL_{\nabla c_i}
+
\sum_{j\in\mathcal E}
|\bar\lambda_j^k|L_{\nabla c_j}
\bigg)
\|d^k\|_2
\\
&\qquad\le
\bigl(
\overline M+\bar L_{\mathcal I}+2L_{\mathcal E}
\bigr)
\|d^k\|_2
=
C_1\|d^k\|_2.
\end{aligned}
\]
Here we used
\[
\bar\lambda_0^k
+
\sum_{i\in\mathcal I}\bar\lambda_i^k
\le1,
\qquad \text{and}\qquad
\sum_{j\in\mathcal E}
|\bar\lambda_j^k|L_{\nabla c_j}
\le
2L_{\mathcal E}.
\]

We next bound the complementarity term.
Since $c_i(x^{k+1}) \le (c_i(\tilde{x}^k))_+$ for every $i\in \mathcal{I}$,
we have
\[
|c_i(x^{k+1})|
\le
(c_i(\tilde{x}^k))_+ + \bigl((c_i(\tilde{x}^k))_+ - c_i(x^{k+1})\bigr),
\qquad i\in \mathcal{I}.
\]
Hence, using $\bar\lambda_i^k \ge 0$, $i\in\mathcal{I}$, we obtain
\begin{align*}
    \sum_{i\in \mathcal{I}} \bar\lambda_i^k |c_i(x^{k+1})| \leq \sum_{i\in \mathcal{I}}\bar\lambda_i^k(c_i(\tilde{x}^k))_+ +\sum_{i\in \mathcal{I}}\bar\lambda_i^k\big((c_i(\tilde{x}^k))_+-c_i(x^{k+1})\big).
\end{align*}
The first sum is bounded by \[
\sum_{i\in\mathcal I}\bar\lambda_i^k(c_i(\tilde{x}^k))_+
\le
\|(c_{\mathcal I}(\tilde x^k))_+\|_2.
\]

For the second sum, the definition of $h_k$ and the Lipschitz continuity of $\nabla f$ give
\begin{align*}
    h_k \geq \nabla f(\tilde{x}^k)^\top d^k +\psi_0(x^{k+1})-\psi_0(\tilde{x}^k) \geq F(x^{k+1})-F(\tilde{x}^k)-(L_{\nabla f}/2)\|d^k\|^2_2.
\end{align*}
Moreover, using the complementary conditions in \eqref{eq:inf-optconds},
\[\bar\lambda_i^k\big(c_i(\tilde{x}^k)+\nabla c_i(\tilde{x}^k)^\top d^k -(c_i(\tilde{x}^k))_+
-h_k\big)=0, \qquad i\in\mathcal{I},\]
together with the Lipschitz continuity of $\nabla c_i$, we obtain
\begin{align*}
    0& \geq \sum_{i\in\mathcal{I}} \bar\lambda_i^k(c_i(x^{k+1})-(c_i(\tilde{x}^k))_+)\\
    & = \sum_{i\in\mathcal{I}} \bar\lambda_i^k(c_i(x^{k+1})-c_i(\tilde{x}^k)-\nabla c_i(\tilde{x}^k)^\top d^k +h_k)\\
    & \geq \sum_{i\in\mathcal{I}} (-\bar\lambda_i^k(L_{\nabla c_i}/2)\|d^k\|^2_2 +\bar\lambda_i^kh_k).
\end{align*}
By \eqref{eq:decrease}, \(h_k\le0\). Since
\(\sum_{i\in\mathcal I}\bar\lambda_i^k\le1\), it follows that $-h_k\sum_{i\in\mathcal I}\bar\lambda_i^k \le-h_k$.
Using the preceding lower bound on \(h_k\), we therefore obtain
\begin{align*}
\sum_{i\in\mathcal{I}}\bar\lambda^{k}_i\bigl((c_i(\tilde{x}^k))_+-c_i(x^{k+1})\bigr)&\leq \bar L_{\mathcal I}\|d^k\|_2^2+ F(\tilde{x}^k)-F(x^{k+1}),
\end{align*}
where we used $\sum_{i\in\mathcal{I}} \bar\lambda_i^k\le 1$, $L_{\nabla f}\leq \bar L_{\mathcal{I}}$, and $L_{\nabla c_i}\leq \bar L_{\mathcal I}$ for all $i\in \mathcal{I}$.
Combining the preceding estimates,
 we obtain 
 \begin{align*}
r^{k+1}_{\mathrm{dist}}
\le \max\bigl\{
C_1\|d^k\|_2, \, 
\bar L_{\mathcal{I}}\|d^k\|_2^2
+F(\tilde{x}^k)-F(x^{k+1})
+\|(c_{\mathcal{I}}(\tilde{x}^k))_+\|_2
\bigr\}.
\end{align*}
Finally, since $\max\{u,v\}\le u+v$, $C_2=\bar L_{\mathcal I}+L_{\mathcal E} \ge\bar L_{\mathcal I}$, and $\|(c_{\mathcal I}(\tilde x^k))_+\|_2\le\|c^+(\tilde x^k)\|_2$, we obtain the second inequality in \eqref{eq:Rk+bound}. The first follows from \Cref{lem:prox-distance-relation}.
\end{proof}

\begin{lemma}
\label{lem:bound_vk}
Let $K\ge0$, and let $a_0,\ldots,a_{K+1}$ and $b_0,\ldots,b_K$ be nonnegative. Suppose that there exist $\rho\in(0,1)$ and $C_{b^2}\ge0$ such that
\begin{equation}\label{eq:sequence-assumptions}
a_{k+1}\le \rho(a_k+b_k),\quad k=0,\ldots,K,\quad
\text{and}\quad
\sum_{k=0}^K b_k^2\le C_{b^2}.
\end{equation}
Then
\begin{equation}
\label{eq:onlyinequality-ak}
\frac{1}{K+1}\sum_{k=0}^K(a_k+b_k)\le\frac{a_0}{(1-\rho)(K+1)}+\frac{1}{1-\rho}\sqrt{\frac{C_{b^2}}{K+1}}.
\end{equation}
\end{lemma}
\begin{proof}
Summing \eqref{eq:sequence-assumptions} over $k=0,\ldots,K$ gives
\[
(1-\rho)\sum_{k=0}^K a_k
\le
a_0+\rho\sum_{k=0}^K b_k.
\]
Combined with the 
Cauchy--Schwarz inequality,
we obtain
\[
\sum_{k=0}^K(a_k+b_k)
\le
\frac{a_0+\sum_{k=0}^K b_k}{1-\rho}
\le
\frac{a_0+\sqrt{(K+1)C_{b^2}}}{1-\rho}.
\]
 Dividing by
$K+1$ proves \eqref{eq:onlyinequality-ak}.
\end{proof}

\subsection{Proof of \Cref{thm:inf-nonasymp}}\label{sec:proof-main-theorem}

\begin{proof}[{Proof of \Cref{thm:inf-nonasymp}}]
Fix $K\ge 0$ and suppose that \Cref{alg:infeasiblefsqp} does not terminate during iterations $0,\ldots,K$.
Define
\[
a_k \coloneqq\|c^+(\tilde x^k)\|_2,\qquad k=0,\ldots,K+1,
\]
and, for \(k=0,\ldots,K\), define
\[
e_{k+1}
\coloneqq
\|c^+(x^{k+1})\|_2,
\qquad
b_k
\coloneqq
C_1\|d^k\|_2
+
C_2\|d^k\|_2^2
+
F(\tilde x^k)-F(x^{k+1}).
\]
By the descent condition \eqref{eq:acceptance-a}, we have $b_k\geq 0$ for all $k \in \{0, \ldots, K\}$.
The residual estimate \eqref{eq:Rk+bound} gives
\begin{equation}
\label{eq:residual-ak-bk}
r^{k+1}_{\mathrm{dist}} \le a_k+b_k.
\end{equation}
Moreover, \eqref{eq:acceptance-b} gives $c_{\mathcal{I}}(x^{k+1})\leq (c_{\mathcal{I}}(\tilde{x}^k))_+$,
componentwise, while
\Cref{lem:basicproperties}\ref{itm:equality-bound} gives
\[
\|c_{\mathcal E}(x^{k+1})\|_2 \le \|c_{\mathcal E}(\tilde x^k)\|_2 + L_{\mathcal E}\|d^k\|_2^2.
\]
Consequently,
\begin{align}
e_{k+1}
&\le
\bigl(
\|(c_{\mathcal I}(\tilde x^k))_+\|_2^2
+
\big(
\|c_{\mathcal E}(\tilde x^k)\|_2
+
L_{\mathcal E}\|d^k\|_2^2
\big)^2
\bigr)^{1/2}
\le
a_k+b_k,
\label{eq:feasibility-ak-bk}
\end{align}
where the last inequality uses $C_2\ge L_{\mathcal E}$ and $F(\tilde x^k)-F(x^{k+1})\ge0$. Hence,
\begin{equation}
\label{eq:max-residual-feasibility}
\max\{r^{k+1}_{\mathrm{dist}},e_{k+1}\} \le a_k+b_k.
\end{equation}

We next derive a recursion for \(a_k\). If the triggering condition
\eqref{eq:triggering-rule} holds, then
\(\tilde x^{k+1}=x^{k+1}\), and therefore
\[
a_{k+1} =e_{k+1} \le \beta r^{k+1}_{\mathrm{prox}}\le \beta r^{k+1}_{\mathrm{dist}} \le \beta(a_k+b_k).
\]
Otherwise, the reduction routine is invoked, and
\Cref{ass:reduction-contraction} gives
\[
a_{k+1} \le \eta_{\mathcal{R}_\theta}e_{k+1} \le \eta_{\mathcal{R}_\theta}(a_k+b_k).
\]
Thus, with
$\rho=\max\{\eta_{\mathcal{R}_\theta},\beta\}$,
\begin{equation}
\label{eq:ak-recursion-main-proof}
a_{k+1} \le \rho(a_k+b_k), \qquad k=0,\ldots,K.
\end{equation}
It remains to bound $\sum_{k=0}^K b_k^2$. By \Cref{lem:basicproperties}\ref{itm:model-decrease} and the objective branch of the max model,
\[
-(L_f+L_{\psi_0})\|d^k\|_2\le h_k\le -M_k\|d^k\|_2^2.
\]
Hence
\[
\|d^k\|_2\le (L_f+L_{\psi_0})/M_k\le (L_f+L_{\psi_0})/\underline M.
\]
Lipschitz continuity of $F$ on $X$ and \eqref{eq:acceptance-a} imply
\[
0\le F(\tilde x^k)-F(x^{k+1})\le (L_f+L_{\psi_0})\|d^k\|_2.
\]
Consequently,
\[
0\le b_k\le\big(C_1+\left(1+C_2/\underline M\right)(L_f+L_{\psi_0})\big)\|d^k\|_2.
\]
Using \eqref{eq:sum-dk} gives
\begin{equation}\label{eq:sum-bk-square}
\sum_{k=0}^K b_k^2\le\big(C_1+\left(1+C_2/\underline M\right)(L_f+L_{\psi_0})\big)^2D_\Theta=Q_\Theta.
\end{equation}
Applying \Cref{lem:bound_vk} to \eqref{eq:ak-recursion-main-proof} and \eqref{eq:sum-bk-square}, and using
\(a_0=\|c^+(x^0)\|_2\), yields
\begin{align*}
    \frac{1}{K+1}\sum_{k=0}^K(a_k+b_k)\leq \frac{\|c^+(x^0)\|_2}{(1-\rho)(K+1)}+\frac{1}{1-\rho}\sqrt{\frac{Q_{\Theta}}{K+1}}.
\end{align*}
By \eqref{eq:max-residual-feasibility},
\[
\begin{aligned}
\min_{0\le k\le K} \max\{r^{k+1}_{\mathrm{dist}},e_{k+1}\} \le \frac{1}{K+1}\sum_{k=0}^K\max\{r^{k+1}_{\mathrm{dist}},e_{k+1}\}\le\frac{1}{K+1}\sum_{k=0}^K(a_k+b_k).
\end{aligned}
\]
Combining this estimate with the preceding averaged bound yields
\[
\min_{0\le k\le K}\max\left\{r^{k+1}_{\mathrm{dist}},\|c^+(x^{k+1})\|_2\right\}\le\frac{\|c^+(x^0)\|_2}{(1-\rho)(K+1)}+\frac{1}{1-\rho}\sqrt{\frac{Q_\Theta}{K+1}},
\]
which proves \eqref{eq:main-nonasymptotic-bound}. The normalized Fritz John proximal residual statement follows from \Cref{lem:prox-distance-relation}.
\end{proof}

\section{Constraint-violation reduction schemes}\label{sec:reduction-step}

In this section, we study two reduction schemes for Step~\ref{step:reduction} of \Cref{alg:infeasiblefsqp}: a proximal-gradient (PG) and a regularized Gauss--Newton (GN) method. Both schemes are applied to the same squared constraint violation function over $X$.

\begin{assumption}[Bounds for the concrete reductions]
\label{ass:constraint-values}
There exist $\overline c_{\mathcal I},\overline c_{\mathcal E}>0$ such that
$|c_i(x)|\le\overline c_{\mathcal I}$ for all $x\in X$, $i\in\mathcal I$, and
$|c_i(x)|\le\overline c_{\mathcal E}$ for all $x\in X$, $i\in\mathcal E$.
\end{assumption}
This additional assumption is used in this section to obtain uniform bounds on the accepted PG and GN backtracking parameters and to relate the stationarity residual at an accepted reduction point to the corresponding PG or GN step measure. It holds automatically for continuous constraint functions on a bounded box $X$. Define
\begin{equation}\label{eq:c^+_max}
c^+_{\max}\coloneqq\big(|\mathcal I|\overline c_{\mathcal I}^{\,2}+|\mathcal E|\overline c_{\mathcal E}^{\,2}\big)^{1/2},
\quad
\text{so that}
\quad \|c^+(x)\|_2\le c^+_{\max}\quad(x\in X).
\end{equation}
The two schemes generate their trial points from different local models.
We first show that their backtracking procedures terminate and that every accepted point lies in $X$ and remains within a distance proportional to the current constraint violation. We then introduce the uniform feasibility--PL condition in \Cref{ass:PL}, which relates a first-order stationarity measure for the squared-violation problem to its function value. Together with the decrease obtained by the reduction step, this condition is used to prove that each scheme reduces the constraint violation by a uniform factor. These results verify all the conditions imposed on the reduction routine in \Cref{ass:reduction-contraction}.

Define the squared constraint-violation function
\begin{equation}\label{eq:prox-gradient-problem}
\upsilon(x)\coloneqq (1/2)\|c^+(x)\|_2^2.
\end{equation}
The following aggregate constant is obtained from
\Cref{ass:standing-regularity}:
\begin{equation}
\label{eq:aggregate-constraint-constants}
\begin{aligned}
L_c
\coloneqq
\left(
|\mathcal I|L_{c,\,\mathcal I}^2
+
|\mathcal E|L_{c,\,\mathcal E}^2
\right)^{1/2}.
\end{aligned}
\end{equation}
By the definitions of \(L_c\) in \eqref{eq:aggregate-constraint-constants} and \(L_{\nabla c}\) in \eqref{eq:aggregate-main-constants}, for all \(x,y\in X\),
\begin{equation}
\label{eq:aggregate-basic-bounds}
\begin{aligned}
\|\nabla c_{\mathcal C}(x)\|_{2}
\le L_c,
\quad \text{and} \quad 
\|\nabla c_{\mathcal C}(x)-\nabla c_{\mathcal C}(y)\|_{2}
\le L_{\nabla c}\|x-y\|_2.
\end{aligned}
\end{equation}
Moreover, $\upsilon$ is continuously differentiable on $X$, and its gradient 
\begin{equation}
\label{eq:gradient-upsilon}
\nabla \upsilon (x) = \nabla c_{\mathcal{C}}(x)c^+(x)
\end{equation}
is Lipschitz continuous on $X$ with constant 
\begin{equation}
\label{eq:L-upsilon}
L_{\nabla\upsilon}
\coloneqq
\overline c_{\mathcal I}\sum_{i\in\mathcal I}L_{\nabla c_i}
+|\mathcal I|L_{c,\,\mathcal I}^2
+\overline c_{\mathcal E}\sum_{i\in\mathcal E}L_{\nabla c_i}
+|\mathcal E|L_{c,\,\mathcal E}^2.
\end{equation}
For $x,y\in X$, define the linearized violation vector
\begin{equation}\label{eq:linearized-violation-vector}
c^+(x;y)
\coloneqq
\begin{bmatrix}
\bigl(c_{\mathcal I}(x)
+\nabla c_{\mathcal I}(x)^\top(y-x)\bigr)_+\\[1mm]
c_{\mathcal E}(x)
+\nabla c_{\mathcal E}(x)^\top(y-x)
\end{bmatrix}.
\end{equation}
Applying the standard Lipschitz-gradient Taylor estimate to each
constraint function, and using the nonexpansiveness of $t\mapsto t_+$, we obtain
\begin{equation}
\label{eq:linearized-violation-error}
\|c^+(y)-c^+(x;y)\|_2
\leq
(L_{\nabla c}/2)\|y-x\|_2^2.
\end{equation}

\subsection{Trial points and common line search}
For $y^0\in X$ and $\alpha>0$, define the PG trial point
\begin{equation}\label{eq:y^+}
        y_{\mathrm {PG}}^+(y^0; \alpha) \coloneqq \mathrm{arg\, min}_{y\in X}\Bigl\{\upsilon(y^0)+\nabla\upsilon(y^0)^\top(y-y^0)+(\alpha/2)\|y-y^0\|_2^2\Bigr\}.
\end{equation}
For the GN trial point, set
\begin{equation}\label{eq:linearized_model}
    \widehat{\upsilon}_{y^0}(y)
\coloneqq (1/2)\| c^+(y^0;y)\|_2^2,
\end{equation}
and define
\begin{equation}
\label{eq:GN-trial-point}
y_{\mathrm{GN}}^+(y^0;\alpha)
\coloneqq
\operatorname*{arg\,min}_{y\in X}
\Bigl\{
\widehat\upsilon_{y^0}(y)
+
(\alpha/2)\|y-y^0\|_2^2
\Bigr\}.
\end{equation}
Both trial points are uniquely defined. 
For $\mathrm{R}\in \{\mathrm{PG},\mathrm{GN}\}$, define
\begin{equation}\label{eq:prox_grad_map}
    \mathcal{H}_{R,\,\alpha}(y^0) \coloneqq \alpha\bigl(y^0-y_{R}^+(y^0;\,\alpha)\bigr).
\end{equation}
The PG and GN trial points are embedded in the same backtracking procedure. The first acceptance condition requires sufficient decrease in the squared constraint violation relative to the corresponding step measure. The second imposes the displacement bound required by \Cref{ass:reduction-contraction}. Algorithm~\ref{alg:prox_gradient} selects the first trial parameter for which both conditions hold.

\begin{algorithm}[!htbp]
\small
\caption{Squared-violation reduction with line search}
\label{alg:prox_gradient}
\begin{enumerate}
    \item[0.] Input $R\in \{\mathrm{PG},\mathrm{GN}\}$, \(y^0\in X\), \(q>1\), \(\alpha_0>0\),
    \(\theta>0\), and \(\sigma\in(0,1)\).

    \item Compute the smallest $\alpha\in \{\alpha_0, q \alpha_0, q^2 \alpha_0,\ldots\}$  for which the corresponding trial point \(y_{R}^+(y^0;\,\alpha)\), defined in \eqref{eq:y^+} and \eqref{eq:GN-trial-point}, respectively,
    satisfies
    \begin{equation}\label{eq:squared-reduction-decrease-test}
    \upsilon\big(y_R^+(y^0;\,\alpha)\big)
    \le
    \upsilon(y^0)
    -
    \frac{\sigma}{2\alpha}
    \|\mathcal H_{R,\,\alpha}(y^0)\|_2^2,
    \end{equation}
    and
    \begin{equation}\label{eq:squared-reduction-displacement-test}
    (1/\alpha)\|\mathcal H_{R,\,\alpha}(y^0)\|_2
    \le
    \theta\|c^+(y^0)\|_2.
    \end{equation}

    \item Set
\[
\widetilde\alpha_{R}\coloneqq \alpha,
\qquad
y_{R}^1\coloneqq y_{R}^+(y^0;\,\widetilde\alpha_{R}),
\]
and return $( \widetilde\alpha_{R},y_{R}^1)$.
\end{enumerate}
\end{algorithm}

For a fixed choice \(R\in\{\mathrm{PG},\mathrm{GN}\}\), the point \(y_R^1\) returned by Algorithm~\ref{alg:prox_gradient} is the candidate output of the corresponding reduction routine. The next two subsections show that the backtracking terminates finitely and that the accepted parameters are uniformly bounded.

\subsection{Proximal-gradient reduction}
For the PG reduction, define 
\begin{equation}
\label{eq:PG-constants}
\begin{aligned}
\underline\alpha_{\mathrm{PG}}
\coloneqq
\max\left\{
\frac{L_c}{\theta},
\frac{L_{\nabla\upsilon}}{1-\sigma}
\right\},
\quad
\overline\alpha_{\mathrm{PG}}
\coloneqq
\max\left\{
\alpha_0,
q\,\underline\alpha_{\mathrm{PG}}
\right\},
\quad
b_{\mathrm{PG}}
\coloneqq
1+\frac{L_{\nabla\upsilon}}{\alpha_0}.
\end{aligned}
\end{equation}
\begin{lemma}[Finite termination and accepted-parameter bound for PG]
\label{lem:PG-backtracking-termination}
Suppose that \Cref{ass:standing-regularity,ass:constraint-values} hold, and let
$y^0\in X$, $q>1$, $\alpha_0>0$, $\theta>0$, and $\sigma\in(0,1)$. Then every $\alpha \geq \underline{\alpha}_{\mathrm{PG}}$ satisfies the acceptance conditions \eqref{eq:squared-reduction-decrease-test}--\eqref{eq:squared-reduction-displacement-test}. Consequently, \Cref{alg:prox_gradient} with $ R=\mathrm{PG}$ terminates finitely, and its accepted parameter $\widetilde \alpha_{\mathrm{PG}}$ satisfies 
 \begin{equation}
     \label{eq:PG-accepted-alpha-bound}
    \alpha_0\leq \widetilde \alpha_{\mathrm{PG}}\leq \overline{\alpha}_{\mathrm{PG}}.
    \end{equation}
    The accepted decrease
\eqref{eq:squared-reduction-decrease-test} and
\eqref{eq:PG-accepted-alpha-bound} imply the uniform estimate
\begin{equation}\label{eq:uniform-upper-reduction_PG}
       \upsilon(y_{\mathrm{PG}}^1) \leq \upsilon(y^0)-\frac{\sigma}{2\overline\alpha_{\mathrm{PG}}}\|\mathcal H_{\mathrm{PG},\widetilde\alpha_{\mathrm{PG}}}(y^0)\|_2^2.
\end{equation}
Moreover, 
\begin{equation}\label{eq:PG-stationarity-transfer}
\operatorname{dist}_2\bigl(-\nabla\upsilon(y_{\mathrm{PG}}^1),N_X(y_{\mathrm{PG}}^1)\bigr)\le b_{\mathrm{PG}}\bigl\|\mathcal H_{\mathrm{PG},\widetilde\alpha_{\mathrm{PG}}}(y^0)\bigr\|_2.
\end{equation}
\end{lemma}
\begin{proof}
    By the definition of $y_{\mathrm {PG}}^+(y^0; \alpha)$ in \eqref{eq:y^+} and comparison with $y^0$, for any $\alpha>0$,
    we have \[
   \nabla\upsilon(y^0)^\top
\bigl(y_{\mathrm{PG}}^+(y^0;\alpha)-y^0\bigr)\leq -(\alpha/2)\|y_{\mathrm {PG}}^+(y^0; \alpha)-y^0\|_2^2.
    \]
    Combining this inequality with the facts that $y^0,\, y_{\mathrm {PG}}^+(y^0; \alpha) \in X$ and Lipschitz continuity of $\nabla \upsilon$, we obtain
    \begin{equation*}
        \begin{aligned}
     \upsilon(y_{\mathrm {PG}}^+(y^0; \alpha))\leq \upsilon(y^0)-\frac{\alpha-L_{\nabla\upsilon}}{2}\|y_{\mathrm {PG}}^+(y^0; \alpha)-y^0\|_2^2= \upsilon(y^0)-\frac{\alpha-L_{\nabla\upsilon}}{2\alpha^2}\|\mathcal{H}_{\mathrm{PG},\,\alpha}(y^0)\|_2^2.
     \end{aligned}
     \end{equation*}
      Thus, the decrease condition \eqref{eq:squared-reduction-decrease-test} holds whenever $\alpha\geq \frac{L_{\nabla\upsilon}}{1-\sigma}$.
      Next, by the standard projection representation of \eqref{eq:y^+} and the nonexpansiveness of the projection,
\[
\|y_{\mathrm {PG}}^+(y^0; \alpha)-y^0\|_2
\le
(1/\alpha)\|\nabla\upsilon(y^0)\|_2.
\]
Under \eqref{eq:aggregate-basic-bounds} and \eqref{eq:gradient-upsilon}, 
\[
\|\nabla\upsilon(y^0)\|_2
\le
L_c\,\|c^+(y^0)\|_2.
\]
Hence,
\[
(1/\alpha)\|\mathcal H_{\mathrm{PG},\,\alpha}(y^0)\|_2
=
\|y_{\mathrm {PG}}^+(y^0; \alpha)-y^0\|_2
\le
(L_c/\alpha)\|c^+(y^0)\|_2.
\]
Therefore \eqref{eq:squared-reduction-displacement-test} holds whenever
$
\alpha\ge \frac{L_c}{\theta}.
$
It follows that every $\alpha\geq \underline{\alpha}_{\mathrm{PG}}$ is accepted.
Since the backtracking procedure generates the sequence \(\alpha_0,q\alpha_0,q^2\alpha_0,\ldots\),  the procedure terminates finitely.
If $\widetilde\alpha_{\mathrm{PG}}>\alpha_0$, then its preceding trial $\widetilde\alpha_{\mathrm{PG}}/q$ was rejected and therefore $\widetilde\alpha_{\mathrm{PG}}<q\underline\alpha_{\mathrm{PG}}$. This proves \eqref{eq:PG-accepted-alpha-bound}.
Moreover, \eqref{eq:uniform-upper-reduction_PG} follows from the accepted decrease test \eqref{eq:squared-reduction-decrease-test} together with $\widetilde\alpha_{\mathrm{PG}}
\leq\overline\alpha_{\mathrm{PG}}$.
Finally, the optimality condition for the PG subproblem gives
\[
\mathcal H_{\mathrm{PG}, \widetilde\alpha_{\mathrm{PG}}}(y^0)-\nabla\upsilon(y^0) \in N_X(y_{\mathrm{PG}}^1).
\]
Therefore,
\[
\begin{aligned}
\operatorname{dist}_2\big(-\nabla\upsilon(y_{\mathrm{PG}}^1), N_X(y_{\mathrm{PG}}^1) \big) &\le
\left\|\mathcal H_{\mathrm{PG},\widetilde\alpha_{\mathrm{PG}}}(y^0) \right\|_2 + L_{\nabla\upsilon} \|y_{\mathrm{PG}}^1-y^0\|_2
\\&=\left(1+\frac{L_{\nabla\upsilon}}{\widetilde\alpha_{\mathrm{PG}}} \right) \left\|\mathcal H_{\mathrm{PG},\, \widetilde\alpha_{\mathrm{PG}}}(y^0)\right\|_2.
\end{aligned}
\]
Then the stationarity estimate \eqref{eq:PG-stationarity-transfer} follows from
\(\widetilde\alpha_{\mathrm{PG}}\ge\alpha_0\).
\end{proof}

\subsection{Regularized Gauss--Newton reduction}\label{sec:GN-reduction}
For the GN reduction, define $\underline\alpha_{\mathrm{GN}}
\coloneqq
\max\left\{
\frac{1}{\theta^2},
\frac{L_{\nabla c}c^+_{\max}}{1-\sigma}
\right\}$, and 
\begin{equation}
\label{eq:GN-constants}
\begin{aligned}
\overline\alpha_{\mathrm{GN}}
\coloneqq
\max\left\{
\alpha_0,
q\,\underline\alpha_{\mathrm{GN}}
\right\},
\quad
b_{\mathrm{GN}}
\coloneqq
1+\frac{L_{\nabla c}c^+_{\max}}{\alpha_0}\left( 1+\frac{\theta L_c}{2}\right).
\end{aligned}
\end{equation}
\begin{lemma}[Finite termination and accepted-parameter bound for GN]\label{lem:GN-backtracking-termination} Suppose that \Cref{ass:standing-regularity,ass:constraint-values} hold, and let $q>1$, $\alpha_0>0$, $\theta>0$, $\sigma\in(0,1)$, and $y^0\in X$. Then every $\alpha\geq\underline\alpha_{\mathrm{GN}}$ satisfies the two acceptance conditions \eqref{eq:squared-reduction-decrease-test}--\eqref{eq:squared-reduction-displacement-test}. Consequently, \Cref{alg:prox_gradient} with $R=\mathrm{GN}$ terminates finitely, and its accepted parameter $\widetilde \alpha_{\mathrm{GN}}$ satisfies 
 \begin{equation}
     \label{eq:GN-accepted-alpha-bound}
    \alpha_0\leq \widetilde \alpha_{\mathrm{GN}}\leq \overline{\alpha}_{\mathrm{GN}}.
    \end{equation}
    The accepted decrease \eqref{eq:squared-reduction-decrease-test} and \eqref{eq:GN-accepted-alpha-bound} imply the uniform estimate
\begin{equation}\label{eq:uniform-upper-reduction_GN}
       \upsilon(y_{\mathrm{GN}}^1) \leq \upsilon(y^0)-\frac{\sigma}{2\overline{\alpha}_{\mathrm{GN}}}\bigl\|\mathcal H_{\mathrm{GN},\widetilde\alpha_{\mathrm{GN}}}(y^0)\bigr\|_2^2.
    \end{equation}
Moreover, 
\begin{equation}
\label{eq:GN-stationarity-transfer}
\operatorname{dist}_2
\bigl(
-\nabla\upsilon(y_{\mathrm{GN}}^1),
N_{X}(y_{\mathrm{GN}}^1)
\bigr)
\le
b_{\mathrm{GN}}
\bigl\|
\mathcal H_{\mathrm{GN},
\widetilde\alpha_{\mathrm{GN}}}(y^0)
\bigr\|_2.
\end{equation}
\end{lemma}
\begin{proof}
Fix $y^0\in X$ and $\alpha>0$. Then $\widehat\upsilon_{y^0}(y^0)=\upsilon(y^0)$, and the optimality of $y_{\mathrm{GN}}^+(y^0;\alpha)$ in \eqref{eq:GN-trial-point} gives 
\begin{equation}
\label{eq:GN-model-comparison}
\widehat\upsilon_{y^0}(y_{\mathrm{GN}}^+(y^0;\alpha))
+
\frac{\alpha}{2}\|y_{\mathrm{GN}}^+(y^0;\alpha)-y^0\|_2^2
\le
\upsilon(y^0).
\end{equation}
Since $\widehat\upsilon_{y^0}(y_{\mathrm{GN}}^+(y^0;\alpha))\geq 0$, 
\begin{equation}
\label{eq:GN-basic-displacement}
\|y_{\mathrm{GN}}^+(y^0;\alpha)-y^0\|_2
\le
(1/\sqrt{\alpha})\|c^+(y^0)\|_2.
\end{equation}
Then acceptance condition \eqref{eq:squared-reduction-displacement-test} holds whenever \(\alpha\geq\theta^{-2}\). 

Moreover, \eqref{eq:GN-model-comparison} implies \[\|c^+(y^0;y_{\mathrm{GN}}^+(y^0;\alpha))\|_2\leq\|c^+(y^0)\|_2\leq c^+_{\max},\] while \(y_{\mathrm{GN}}^+(y^0;\alpha)\in X\) gives
\[
\|c^+(y_{\mathrm{GN}}^+(y^0;\alpha))\|_2
\le
c^+_{\max}.
\] Thus, by \eqref{eq:linearized-violation-error}, the difference-of-squares identity, and the reverse triangle inequality,
\begin{equation}\label{eq:GN-actual-model}
\upsilon(y_{\mathrm{GN}}^+(y^0;\alpha))-\widehat\upsilon_{y^0}(y_{\mathrm{GN}}^+(y^0;\alpha))
\leq
\frac{L_{\nabla c}\,c^+_{\max}}{2}
\|y_{\mathrm{GN}}^+(y^0;\alpha)-y^0\|_2^2.
\end{equation}

 Combining
\eqref{eq:GN-model-comparison} and \eqref{eq:GN-actual-model} gives
\begin{align*}
\upsilon(y_{\mathrm{GN}}^+(y^0;\alpha))
&\leq
\upsilon(y^0)
-
\bigg(\frac{\alpha}{2}-\frac{L_{\nabla c}\,c^+_{\max}}{2}\bigg)
\|y_{\mathrm{GN}}^+(y^0;\alpha)-y^0\|_2^2
\\
&=
\upsilon(y^0)
-
\frac{\alpha-L_{\nabla c}\,c^+_{\max}}{2\alpha^2}
\|\mathcal H_{\mathrm{GN},\,\alpha}(y^0)\|_2^2.
\end{align*}
The decrease acceptance condition \eqref{eq:squared-reduction-decrease-test} therefore holds whenever
\(\alpha\geq L_{\nabla c}\,c^+_{\max}/(1-\sigma)\).  Thus every \(\alpha\geq\underline\alpha_{\mathrm{GN}}\) is accepted.  The geometric
backtracking sequence gives \eqref{eq:GN-accepted-alpha-bound}, and the
accepted decrease test \eqref{eq:squared-reduction-decrease-test} together with $\widetilde\alpha_{\mathrm{GN}}
\leq\overline\alpha_{\mathrm{GN}}$ give \eqref{eq:uniform-upper-reduction_GN}.

Finally, using \eqref{eq:gradient-upsilon},
\(\nabla\widehat\upsilon_{y^0}(y_{\mathrm{GN}}^1)
=
\nabla c_{\mathcal C}(y^0)c^+(y^0;y_{\mathrm{GN}}^1)\),
and adding and subtracting
\(\nabla c_{\mathcal C}(y^0)c^+(y_{\mathrm{GN}}^1)\), we obtain, by \eqref{eq:c^+_max},
\eqref{eq:aggregate-basic-bounds} and
\eqref{eq:linearized-violation-error}
\[
\begin{aligned}
\|\nabla\upsilon(y_{\mathrm{GN}}^1)
-
\nabla\widehat\upsilon_{y^0}(y_{\mathrm{GN}}^1)\|_2
&\leq \|\bigl(
\nabla c_{\mathcal C}(y_{\mathrm{GN}}^1)
-\nabla c_{\mathcal C}(y^0)
\bigr)c^+(y_{\mathrm{GN}}^1)\|_2
\\&\quad+\|\nabla c_{\mathcal C}(y^0)
\bigl(
c^+(y_{\mathrm{GN}}^1)
-
c^+(y^0;y_{\mathrm{GN}}^1)
\bigr)\|_2\\
& \leq L_{\nabla c}c^{+}_{\max}\|y^1_{\mathrm{GN}}-y^0\|_2 + \frac{L_cL_{\nabla c}}{2}\|y^1_{\mathrm{GN}}-y^0\|_2^2.
\end{aligned}
\]
Moreover, since $y^1_{\mathrm{GN}}$ satisfies $\|y^1_{\mathrm{GN}}-y^0\|_2\leq \theta \|c^+(y^0)\|_2\leq \theta c^+_{\max}$, we have
\[
\|\nabla\upsilon(y_{\mathrm{GN}}^1)
-
\nabla\widehat\upsilon_{y^0}(y_{\mathrm{GN}}^1)\|_2 \leq L_{\nabla c}c^+_{\max}
\left(
1+\frac{\theta L_c}{2}
\right)\|y_{\mathrm{GN}}^1-y^0\|_2.
\]
The GN subproblem optimality condition gives
\[
\mathcal H_{\mathrm{GN},
\widetilde\alpha_{\mathrm{GN}}}(y^0)
-
\nabla\widehat\upsilon_{y^0}(y_{\mathrm{GN}}^1)
\in
N_X(y_{\mathrm{GN}}^1).
\]
Hence we have
\[
\begin{aligned}
\operatorname{dist}_2\left(-\nabla\upsilon(y_{\mathrm{GN}}^1), N_{X}(y_{\mathrm{GN}}^1)\right)&\le\left\|\mathcal H_{\mathrm{GN},\widetilde\alpha_{\mathrm{GN}}}(y^0)\right\|_2+L_{\nabla c}c^+_{\max}\left(1+\frac{\theta L_c}{2}\right)\|y_{\mathrm{GN}}^1-y^0\|_2\\&\leq b_{\mathrm{GN}}\left\|\mathcal H_{\mathrm{GN}, \widetilde\alpha_{\mathrm{GN}}}(y^0)\right\|_2, 
\end{aligned}
\]
where the last inequality follows from \(\widetilde\alpha_{\mathrm{GN}}\ge\alpha_0\).
\end{proof}

\subsection{Uniform feasibility--PL condition}
The preceding PG and GN analyses, under \Cref{ass:constraint-values}, provide sufficient decrease in the squared constraint violation and stationarity estimates at the accepted reduction points. The following PL condition links these stationarity estimates to the squared constraint violation. Combined with the accepted decrease, it yields a uniform contraction of the constraint violation.

\label{sec:uniform-feasibility-PL}
\begin{assumption}[Uniform PL condition]\label{ass:PL}
There exists a constant $\mu>0$ such that 
\[
\operatorname{dist}_2^2\!\bigl(-\nabla \upsilon(x),\,N_{X}(x)\bigr) \ge 2\mu\upsilon(x)\quad \text{for all} \quad  x\in X.
\]
\end{assumption}

\begin{remark}[Verifying the feasibility--PL condition]
\label{rem:verifying-PL}
Recall that \(S\) denotes the feasible set of \eqref{eq:infeasible-fsqp}. Suppose that \(\upsilon\) is convex on \(X\) and satisfies the quadratic-growth condition
\[
\upsilon(x) \ge2\mu\,\operatorname{dist}_2^2(x,S)\quad \text{for all} \quad  
x\in X.
\]
Then \Cref{ass:PL} follows from the standard relationship between quadratic growth and subdifferential error bounds for convex functions, applied to \(\upsilon+\delta_X\); see \cite[Theorem~3.3]{Drusvyatskiy2018}. Related relationships between quadratic growth and PL inequalities are discussed in \cite[Theorem~3.1]{LiaoDingZheng2024}. An example illustrating that the feasibility--PL condition can hold despite failure of LICQ and MFCQ is provided in \Cref{app:verifying-PL}.
\end{remark}

\begin{theorem}[Uniform admissibility of PG and GN]
\label{thm:squared-reduction-admissibility}
Suppose that \Cref{ass:standing-regularity,ass:constraint-values,ass:PL} hold. Fix \(q>1\), \(\alpha_0>0\), \(\theta>0\), and \(\sigma\in(0,1)\). For each fixed choice \(R\in\{\mathrm{PG},\mathrm{GN}\}\), let \(\mathcal R_\theta:X\to X\) denote the output mapping of \Cref{alg:prox_gradient}, that is,
$
\mathcal R_\theta(y^0)\coloneqq y_R^1
$.
Then \(\mathcal R_\theta\) satisfies \Cref{ass:reduction-contraction} with
\begin{equation}
\label{eq:squared-reduction-contraction-factor}
\eta_{\mathcal R_\theta}=\eta_R\coloneqq\left(1+\frac{\sigma\mu}{\overline\alpha_R b_R^2}\right)^{-1/2}\in(0,1),
\end{equation}
where \(\mu\) is the constant in \Cref{ass:PL}, and \(\overline\alpha_R\) and \(b_R\) are the corresponding constants in \eqref{eq:PG-constants} or \eqref{eq:GN-constants}.
\end{theorem}

\begin{proof}
Fix \(R\in\{\mathrm{PG},\mathrm{GN}\}\) and \(y^0\in X\). By \Cref{lem:PG-backtracking-termination} when \(R=\mathrm{PG}\) and \Cref{lem:GN-backtracking-termination} when \(R=\mathrm{GN}\), \Cref{alg:prox_gradient} terminates finitely. Its returned point \(y_R^1\) lies in \(X\), and \eqref{eq:squared-reduction-displacement-test} together with \eqref{eq:prox_grad_map} gives \(\|y_R^1-y^0\|_2 \le \theta\|c^+(y^0)\|_2\). 

It remains to establish the contraction bound.
By \eqref{eq:uniform-upper-reduction_PG} when \(R=\mathrm{PG}\) and \eqref{eq:uniform-upper-reduction_GN} when \(R=\mathrm{GN}\),
\[
\upsilon(y_R^1)
\le
\upsilon(y^0)
-
\frac{\sigma}{2\overline\alpha_R}
\|\mathcal H_{R,\widetilde\alpha_R}(y^0)\|_2^2.
\]
On the other hand, \Cref{ass:PL}, together with
\eqref{eq:PG-stationarity-transfer} for \(R=\mathrm{PG}\) and
\eqref{eq:GN-stationarity-transfer} for \(R=\mathrm{GN}\), gives
\[
\begin{aligned}
2\mu\upsilon(y_R^1)\le\operatorname{dist}_2^2\bigl(-\nabla\upsilon(y_R^1),N_X(y_R^1)\bigr)\le b_R^2 \bigl\|\mathcal H_{R,\widetilde\alpha_R}(y^0)\bigr\|_2^2.
\end{aligned}
\]
Combining these inequalities yields
$
\upsilon(y_R^1)
\le
\eta_R^2
\upsilon(y^0)
$.
Since
\((1/2)\|c^+(x)\|_2^2=\upsilon(x)\), we have
\( \|c^+(y_R^1)\|_2\le\eta_R\|c^+(y^0)\|_2\).
Since \(\eta_R\in(0,1)\) is independent of \(y^0\), the output
mapping satisfies \Cref{ass:reduction-contraction}.
\end{proof}
If, in addition, \Cref{ass:obj-lower-bound} holds, then \Cref{thm:inf-nonasymp} applies to \Cref{alg:infeasiblefsqp} equipped with either reduction scheme.

\section{Numerical results}\label{sec:Numerical_implementation}

We present numerical experiments using a Julia \cite{Bezanson2017Julia} implementation of \Cref{alg:infeasiblefsqp}. We examine multiplier behavior on a degenerate drinking-water network problem, compare reduction schemes and a curvature-enhanced variant on a subset of the \textsc{CUTEst} test problems \cite{GouldOrbanToint2015CUTEst}, and compare \ISMFD{} with \LCPG{} \cite{Boob2025} on \textsc{QCQP} instances. All experiments were run on the Georgia Tech PACE Phoenix cluster~\cite{PACE}. The \textsc{CUTEst} experiments used four allocated CPU cores and $8$~GB of memory per problem instance, while the \textsc{QCQP} comparison used one core and $16$~GB per instance. We solve the subproblems in \eqref{eq:subproblem'} using \Clarabel{} \cite{GoulartChen2026Clarabel}, with a maximum of $200$ iterations, absolute and relative duality-gap tolerances of $10^{-10}$, and the default feasibility tolerance $10^{-8}$. Throughout the numerical experiments, we use \(\tau_0=1\), \(\gamma_\tau=2\), \(\underline M=0.1\), and \(\beta=10^{-3}\).

\paragraph{Warm-started regularization search.}
In the numerical experiments, we use a warm-started variant of the first-acceptable doubling search in Step~\ref{step:main}.  We set \(\widehat M_0=2\underline M\), start iteration \(k\) from
\[
M_{k,0} = \max\{\underline M,\widehat M_k\},
\]
and choose the smallest $j\geq 0$ for which 
$M_{k,j}=2^jM_{k,0}$ satisfies \eqref{eq:acceptance-a}--\eqref{eq:acceptance-c}. The accepted value is denoted by \(M_k\), with corresponding step \(d^k\) and point \(x^{k+1}=\tilde x^k+d^k\).
To initialize the next regularization search, we use the agreement between the predicted and actual objective decrease. Define
\[
\operatorname{pred}_k
\coloneqq
-\nabla f(\tilde x^k)^\top d^k
+\psi_0(\tilde x^k)-\psi_0(x^{k+1}),
\qquad
\zeta_k
\coloneqq
\frac{F(\tilde x^k)-F(x^{k+1})}
     {\operatorname{pred}_k},
\]
where \(\zeta_k\) is evaluated only when
\(\operatorname{pred}_k>0\).  The next carry-over value is
\[
\widehat M_{k+1}
=
\begin{cases}
2M_k,
&
\operatorname{pred}_k>0
\ \text{and}\ 
\zeta_k<0.1,
\\[1mm]
M_k/2,
&
\operatorname{pred}_k>0,\ 
\zeta_k\ge0.9,
\ \text{and \(\mathcal{R}_\theta\) is not invoked at iteration \(k\)},
\\[1mm]
M_k,
&
\text{otherwise}.
\end{cases}
\]
For the exact-solve analysis of the scalar-regularized main subproblem, the warm-started search retains the uniform bound \(\underline M\le M_k\le\overline M\).
The proof is given in \Cref{rem:warm-start-bound}. The numerical implementation instead uses finite-precision \Clarabel{} solves.

\begin{remark}[Interpretation of the warm-start update]
\label{rem:M-warm-start-interpretation}
The update follows the model-agreement principle in trust-region methods \cite[Chapter~4]{NocedalWright2006}. When the predicted decrease is positive, poor agreement (\(\zeta_k<0.1\)) motivates doubling the next initial regularization to discourage overly large steps, even though the current step has passed all three acceptance tests. By contrast, when \(\zeta_k\geq0.9\) and no reduction is invoked, we halve the next initial regularization to allow larger trial steps. Retaining strong regularization in this case may unnecessarily slow progress along directions of low curvature. For intermediate agreement (\(0.1\leq\zeta_k<0.9\)), we retain the current value. These adjustments affect only the starting value of the next search, whose accepted step must still satisfy all three acceptance conditions.
\end{remark}

\paragraph{Stopping criterion.}
If $\bar\lambda_0^k>10^{-12}$, 
we set $\widehat\lambda^k=\bar\lambda^k/\bar\lambda_0^k$.
Otherwise, we formally set $r_{\mathrm{prox}}
(
x^{k+1},
1,
\widehat\lambda^k
) = +\infty$.
We terminate \Cref{alg:infeasiblefsqp} if
\[
\|c^+(x^{k+1})\|_2
\le
\varepsilon_{\mathrm{feas}},
\qquad
r_{\mathrm{prox}}
(
x^{k+1},
1,
\widehat\lambda^k
)
\le
\varepsilon_{\mathrm{stat}},
\]
where $\varepsilon_{\mathrm{feas}}=10^{-5}$ and $\varepsilon_{\mathrm{stat}}=10^{-4}$. 
The reduction trigger remains unchanged from \Cref{alg:infeasiblefsqp} and uses \(r_{\mathrm{prox}}(x^{k+1},\bar\lambda_0^k,\bar\lambda^k)\).

\paragraph{Reduction schemes.}
We compare the PG and GN reductions analyzed in \Cref{sec:reduction-step} with an experimental prox-linear reduction, denoted by ProxLin.
 Given \(y^0\in X\) and \(\alpha>0\), its trial point is
\[
y_{\mathrm{ProxLin}}^+(y^0;\alpha)
\coloneqq
\operatorname*{arg\, min}_{y\in X}
\left\{
\|c^+(y^0;y)\|_2
+
\frac{\alpha}{2}\|y-y^0\|_2^2
\right\}.
\]
A trial \(y^+=y_{\mathrm{ProxLin}}^+(y^0;\alpha)\) is accepted if
$$ \|c^+(y^+)\|_2 \le \|c^+(y^0)\|_2 - \frac{\sigma\alpha}{2}\|y^+-y^0\|_2^2 $$ and
\[
\|y^+-y^0\|_2
\le
\theta\|c^+(y^0)\|_2.
\]
If either test fails, \(\alpha\) is multiplied by \(q\). ProxLin is not covered by the PG/GN complexity result. 
For all three reduction routines, we use \(\theta=500\), \(q=2\), \(\alpha_0=10^{-3}\), and \(\sigma=0.5\).

\subsection{Comparison with LCPG on QCQP instances}\label{sec:QCQP}

Following \cite[Section~8.1]{Boob2025}, we consider the QCQP family
\begin{equation}
\label{eq:QCQP-family}
\begin{aligned}
\min_{x\in\mathbb R^n}\quad
&
(1/2) x^\top Q_0x + b_0^\top x + \|x\|_1
\\
\text{s.t.}\quad
&
(1/2) x^\top Q_ix + b_i^\top x -10 \le0,
\quad
i=1,\ldots,9,
\quad \|x\|_2^2-20\le 0,
\end{aligned}
\end{equation}
where $Q_i=V_iD_iV_i^\top$ with $V_i\in\mathbb R^{n\times n}$ a sparse random matrix of density $0.01$ and nonzeros uniform on $[0,1]$, and $D_i \in\mathbb R^{n\times n}$ diagonal with entries uniform on $[0, 100]$. Moreover, $b_i=10e+v_i$, where \(e\in\mathbb R^n\) is the vector of all ones, and \(v_i\sim\mathcal N(0,I_n)\). For nonconvex instances, we instead set \(Q_i=V_i D_i V_i^\top-10I_n\).

\paragraph{Experimental setup.} For each considered $n \in \mathbb{N}$, we generate five random instances. Both methods start from \(x^0=0\) and are run for at most \(1000\) iterations.  The reported means and medians include runs that reach the iteration limit. \ISMFD{} uses PG reduction, while
\LCPG{} is implemented in Julia following \cite{Boob2025}. For a common stopping criterion, we use \(r_{\mathrm{prox}}(x,1,\widehat\lambda_{\mathcal I})\) for both methods. 
For \ISMFD{},
\(\widehat\lambda_{\mathcal I}
=\bar\lambda_{\mathcal I}/\bar\lambda_0\),
whereas for \LCPG{},
\(\widehat\lambda_{\mathcal I}\) is given by its returned inequality
multipliers.

\paragraph{Convex and nonconvex instances.}
As shown in \Cref{tab:qcqp}, \ISMFD{} succeeds on all \(25\) convex instances, compared with $10$ for \LCPG{}. On the nonconvex family, \ISMFD{} succeeds on \(24\) of the \(25\) instances, compared with $8$ for \LCPG{}. In both families, the reported mean objective values of the two methods agree to the displayed precision. Across all reported dimensions, \ISMFD{} has lower median iteration counts and wall-clock times.

\begin{table}[t]
\centering
\footnotesize
\setlength{\tabcolsep}{4pt}
\caption{QCQP comparison using five instances per dimension. Objectives are
means, while iterations, times, and residuals are medians. Successes count
instances satisfying both stationarity and feasibility tolerances within the
$1000$-iteration cap. An asterisk marks a residual reported at the iteration cap.}
\label{tab:qcqp}
\begin{tabular}{rl rrrr c rrrr c}
\toprule
 & & \multicolumn{5}{c}{convex} & \multicolumn{5}{c}{nonconvex} \\
\cmidrule(lr){3-7}\cmidrule(lr){8-12}
$n$ & method & obj & iters & time & $r_{\mathrm{prox}}$ & (succ.)
              & obj & iters & time & $r_{\mathrm{prox}}$ & (succ.) \\
\midrule
  \multirow{2}{*}{500} & \ISMFD{} & $-166.64$ & 114 & 3\,s & $8.7\times10^{-5}$ & (5/5) & $-201.24$ & 485 & 10\,s & $8.2\times10^{-5}$ & (4/5) \\
   & \LCPG{} & $-166.64$ & 812 & 149\,s & $9.9\times10^{-5}$ & (3/5) & $-201.24$ & 967 & 194\,s & $9.8\times10^{-5}$ & (3/5) \\
\addlinespace
  \multirow{2}{*}{1000} & \ISMFD{} & $-85.22$ & 109 & 9\,s & $7.6\times10^{-5}$ & (5/5) & $-91.05$ & 165 & 12\,s & $9.5\times10^{-5}$ & (5/5) \\
   & \LCPG{} & $-85.22$ & 575 & 261\,s & $9.7\times10^{-5}$ & (5/5) & $-91.05$ & 727 & 335\,s & $9.8\times10^{-5}$ & (5/5) \\
\addlinespace
  \multirow{2}{*}{2000} & \ISMFD{} & $-44.07$ & 102 & 42\,s & $7.7\times10^{-5}$ & (5/5) & $-45.33$ & 94 & 33\,s & $6.6\times10^{-5}$ & (5/5) \\
   & \LCPG{} & $-44.07$ & 1000 & 983\,s & $1.2\times10^{-4}$\rlap{$^{*}$} & (2/5) & $-45.33$ & 1000 & 988\,s & $3.6\times10^{-4}$\rlap{$^{*}$} & (0/5) \\
\addlinespace
  \multirow{2}{*}{3000} & \ISMFD{} & $-30.02$ & 135 & 122\,s & $9.3\times10^{-5}$ & (5/5) & $-30.58$ & 141 & 125\,s & $8.6\times10^{-5}$ & (5/5) \\
   & \LCPG{} & $-30.02$ & 1000 & 1695\,s & $3.0\times10^{-3}$\rlap{$^{*}$} & (0/5) & $-30.58$ & 1000 & 1622\,s & $4.5\times10^{-3}$\rlap{$^{*}$} & (0/5) \\
\addlinespace
  \multirow{2}{*}{4000} & \ISMFD{} & $-22.68$ & 125 & 240\,s & $5.6\times10^{-5}$ & (5/5) & $-23.02$ & 126 & 241\,s & $8.3\times10^{-5}$ & (5/5) \\
   & \LCPG{} & $-22.68$ & 1000 & 2316\,s & $2.1\times10^{-2}$\rlap{$^{*}$} & (0/5) & $-23.02$ & 1000 & 2168\,s & $2.7\times10^{-2}$\rlap{$^{*}$} & (0/5) \\
\bottomrule
\end{tabular}
\end{table}

\subsection{Multiplier behavior under degeneracy}
\label{sec:degenerate-examples}
We consider the drinking-water network problem \cite[Section 4.2]{ HaeserHinderYe2021}, with decision variables \(y,z\in\mathbb R^3\). We write \(x\coloneqq(y,z)\in\mathbb R^6\) for the full decision vector.

\begin{equation}\label{eq:drink}
\begin{aligned}
\min_{y,\,z\,\ge\,0}\; z_1 \quad \text{s.t.}\quad
& y_1+y_2=1, && y_1+y_3=0.5, && y_2=0.5,\\
& z_1-z_2=y_1^{1.8}, && z_1-z_3=y_2^{1.8}, && y_3^{1.8}=z_2-z_3.
\end{aligned}
\end{equation}
Problem \eqref{eq:drink} has a unique global solution with $z_1^\star=2^{-1.8}$, at which the KKT conditions hold but MFCQ fails. Consequently, the KKT multiplier set is unbounded. In \cite[Section~4.2]{HaeserHinderYe2021}, this example is used to illustrate large dual multiplier estimates returned by \IPOPT{}.

All methods are initialized at $x^0=10^{-2}e$. \IPOPT{} \cite{WachterBiegler2006} is run with \texttt{tol=1e-4} and a $5000$-iteration limit.
For both methods, the last column of \Cref{tab:drink} reports KKT proximal residual $r_{\mathrm{prox}}(x,1,\widehat\lambda_{\mathcal E})$ at their respective returned primal-dual pairs.  At the terminal \ISMFD{} iteration \(k\), we set $x=x^{k+1}$, and $\widehat\lambda_{\mathcal E} = \bar\lambda_{\mathcal E}^k/\bar\lambda_0^k$.  For \IPOPT{}, \(x\) and \(\widehat\lambda_{\mathcal E}\) are its returned primal point and equality multiplier estimate, respectively. At the final iterate of each \ISMFD{} run, $\bar\lambda_0^k\approx0.365$ for ProxLin and $\bar\lambda_0^k\approx0.331$ for both PG and GN.

As shown in \Cref{tab:drink}, the default \IPOPT{} run reports successful termination, while the no-bound-relaxation run exits after finding a feasible point. Both runs return very large equality-multiplier estimates and large KKT proximal residuals. In contrast, all three \ISMFD{} variants return much smaller multiplier norms and KKT residuals, together
with smaller objective errors.

\begin{table}[t]
\centering
\caption{Diagnostics for the water-network problem \eqref{eq:drink}.}
\label{tab:drink}
\resizebox{\linewidth}{!}{%
\begin{tabular}{llccccc}
\toprule
Method
& Status
& Iter.
& \(\|\widehat\lambda_{\mathcal E}\|_\infty\)
& \(\lvert z_1-z_1^\star\rvert\)
& \(\|c^+(x)\|_2\)
& \(r_{\mathrm{prox}}(x,1,\widehat\lambda_{\mathcal E})\)
\\
\midrule
\ISMFD{} (ProxLin)
& Stopping criterion met
& 3
& \(5.9\times10^{-1}\)
& \(3.0\times10^{-11}\)
& \(1.8\times10^{-10}\)
& \(1.6\times10^{-7}\)
\\
\ISMFD{} (PG)
& Stopping criterion met
& 47
& \(6.2\times10^{-1}\)
& \(1.2\times10^{-6}\)
& \(9.2\times10^{-6}\)
& \(3.8\times10^{-6}\)
\\
\ISMFD{} (GN)
& Stopping criterion met
& 4
& \(6.2\times10^{-1}\)
& \(7.4\times10^{-7}\)
& \(5.9\times10^{-6}\)
& \(3.2\times10^{-5}\)
\\
\textsc{Ipopt} (no bound relax.)
& Feasible-point exit
& 6
& \(3.0\times10^{14}\)
& \(4.1\times10^{-4}\)
& \(7.5\times10^{-7}\)
& \(2.2\times10^{8}\)
\\
\textsc{Ipopt} (default)
& Solved Succeeded
& 19
& \(5.0\times10^{11}\)
& \(1.2\times10^{-4}\)
& \(3.8\times10^{-10}\)
& \(0.29\)
\\
\bottomrule
\end{tabular}%
}
\smallskip
\parbox{0.96\linewidth}{\footnotesize
The no-bound-relaxation run uses
\texttt{bound\_relax\_factor\(=0\)}.
Its ``Feasible-point exit'' is \textsc{Ipopt}'s
\texttt{Feasible\_Point\_Found} status.}
\end{table}

\subsection{Benchmark on the \textsc{CUTEst} collection}
\label{sec:cutest}
We consider all twice continuously differentiable \textsc{CUTEst} problems with an objective, at least one general constraint, and \(1\le n\le200\). We exclude four instances with nonfinite function or derivative values at the initial point and four with more equality constraints than non-fixed variables. The resulting \(344\) problems are solved from their \textsc{CUTEst} initial points with a \(50000\)-iteration cap and a one-hour time limit.

\paragraph{Curvature-enhanced variant.}
We compare the scalar-regularized formulation \eqref{eq:subproblem'}
with an experimental curvature-enhanced variant.
We set \(B_0=0\) and, for \(k\ge1\), define
\[
B_k
=
\lambda_0^{k-1}\nabla^2 f(\tilde x^k)
+
\sum_{i\in\mathcal I}\lambda_i^{k-1}\nabla^2 c_i(\tilde x^k)
+
\sum_{j\in\mathcal E}\lambda_j^{k-1}\nabla^2 c_j(\tilde x^k),
\]
where the multipliers are the \Clarabel{} dual multipliers from the preceding accepted main step, without additional normalization involving the equality multipliers. The regularization search uses the same carry-over warm-start rule for \(\widehat M_k\) as the scalar formulation and starts at
\[
M_{k,0} = \max\bigl\{\underline M,\,\widehat M_k,\,10^{-8}-\lambda_{\min}(B_k)\bigr\}.
\]
For each trial value \(M_k\), we set \(Q_k\coloneqq B_k+M_kI\succeq10^{-8}I\) and solve \eqref{eq:subproblem'} with \((M_k/2)\|d\|_2^2\) replaced by \((1/2)d^\top Q_kd\). The resulting step \(d^k\) is checked against \eqref{eq:acceptance-a}--\eqref{eq:acceptance-c}, with \(M_k\|d^k\|_2^2\) replaced by \((d^k)^\top Q_kd^k\) in \eqref{eq:acceptance-a} and \eqref{eq:acceptance-c}, while \eqref{eq:acceptance-b} is unchanged. If any test fails, we double \(M_k\), update \(Q_k\), and solve the subproblem again. Otherwise, the current \(M_k\) and \(d^k\) are accepted.
This experimental variant is not covered by the scalar complexity analysis or its warm-start parameter bound. In particular, no uniform upper bound is established for the spectral shift \(10^{-8}-\lambda_{\min}(B_k)\).

\begin{table}[t]
\centering
\caption{\textsc{CUTEst} benchmark.}
\label{tab:cutest}
\begin{tabular}{lccc}
\toprule
Reduction
& Problems
& Scalar formulation
& Curvature variant
\\
\midrule
PG
& 344
& 259 (75.3\%)
& 265 (77.0\%)
\\
GN
& 344
& 265 (77.0\%)
& 276 (80.2\%)
\\
ProxLin
& 344
& 282 (82.0\%)
& 294 (85.5\%)
\\
\bottomrule
\end{tabular}
\end{table}

\paragraph{Results.} \Cref{tab:cutest} reports the number and percentage of successful runs. For both formulations, the solved count increases from PG, GN, to ProxLin, and the curvature-enhanced variant improves the count for each reduction scheme. All unsuccessful runs reach the $50000$-iteration cap or the one-hour time limit.

In PG, GN, and ProxLin order, the iteration-limit/time-limit counts are $84/1$, $43/36$, and $51/11$ for the scalar formulation, and $78/1$, $31/37$, and $41/9$ for the curvature-enhanced variant.

\begin{figure}[t]
\centering
\includegraphics[width=\textwidth]{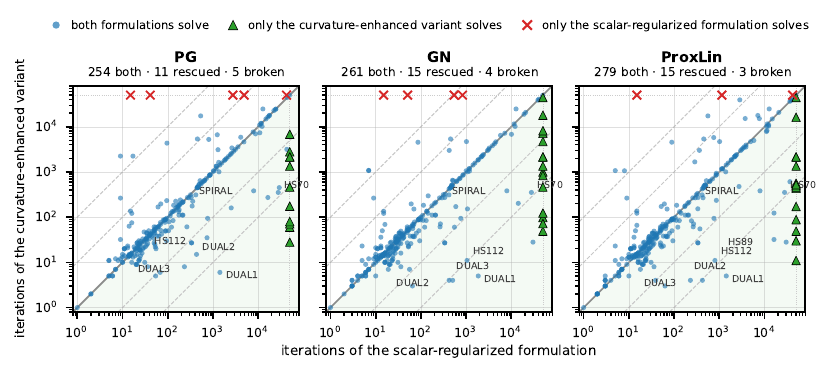}
\caption{Iteration counts for the scalar-regularized formulation and the curvature-enhanced variant on \textsc{CUTEst} (log scales). Unsuccessful runs are assigned $50000$ iterations for plotting purposes. Points below the diagonal favor the curvature-enhanced variant, and dashed lines indicate \(10\times\) and \(100\times\) differences. Triangles and crosses denote problems solved only by the curvature-enhanced and scalar-regularized formulations, respectively.}
\label{fig:so-scatter}
\end{figure}

\Cref{fig:so-scatter} compares the iteration counts of two formulations. Most problems lie near the diagonal, indicating similar iteration counts on many problems,
while the curvature-enhanced variant gives substantial reductions on a small subset. Together with the higher solved counts in \Cref{tab:cutest}, this suggests that its main benefit is improved robustness, with occasional substantial reductions in iteration count.

\section{Discussion}
\label{sec:conclusion}
We introduced \ISMFD{}, an infeasible-start feasible-direction framework combining a regularized max-type main step with a separate feasibility-reduction routine. Under the stated assumptions, including uniform contraction and displacement bounds for the reduction routine, we established an \(O(\varepsilon^{-2})\) outer-iteration bound for generating an iterate whose normalized Fritz John proximal residual and constraint violation are both at most $\varepsilon$. The analysis uses a stronger best-iterate distance-residual estimate and a uniform bound on the selected regularization parameters. We impose neither LICQ nor MFCQ, but the reduction property entails a global feasibility error bound. We verified this property for PG and GN under the feasibility--PL condition and additional constraint-value bounds.

The numerical results illustrate the flexibility of the framework. On the degenerate drinking-water network problem, the \ISMFD{} variants returned small KKT proximal residuals and substantially smaller
multiplier norms than the reported \IPOPT{} runs. On \textsc{CUTEst}, the curvature-enhanced variant increased the solved count and substantially reduced iteration counts on some instances. On the tested convex and nonconvex QCQP families, \ISMFD{} satisfied the common stationarity and feasibility stopping criteria more quickly and reliably than \LCPG{}.

The main step controls possible increases in equality violation rather than enforcing its monotone decrease. Future work includes safeguarded equality corrections that preserve sufficient objective decrease and inequality-violation control, weaker feasibility assumptions verified on invariant regions, and controlled inexactness in the main and reduction subproblems. Further directions include extending the box constraint to broader closed convex sets and analyzing multiplier selection when classical constraint qualifications fail.

\section*{Acknowledgments} 
We thank Professor Digvijay Boob for sharing his implementation of \LCPG{}.

\paragraph{Use of generative AI.}
Generative AI tools were used in a limited supporting role for language editing, presentation, and selected coding and debugging tasks. The research questions, methodology, mathematical analysis, experimental design, and interpretation of the results were developed by the authors.

\appendix
\section{Compatibility of the feasibility--PL condition with constraint-qualification failure}
\label[appendix]{app:verifying-PL} 
The following example illustrates that the feasibility--PL condition in \Cref{ass:PL}, used in the PG/GN contraction analysis, can hold even when LICQ and MFCQ fail. Recall $\upsilon$ from \eqref{eq:prox-gradient-problem} and the feasible set $S$ of \eqref{eq:infeasible-fsqp}.

\begin{remark}[Compatibility with constraint-qualification failure]
\label{rem:PL-without-CQ}
Consider \(X=[-1,1]\) with the redundant equality constraints
\[
c_1(x)=x, \qquad c_2(x)=2x.
\]
The feasible set is \(S=\{0\}\). At the feasible point, the two equality gradients are linearly dependent, so MFCQ, and consequently LICQ, fail.
Nevertheless,
\[
\operatorname{dist}_2(x,S) = |x| \le \frac{1}{\sqrt5} \left\|
\begin{bmatrix}
c_1(x)\\
c_2(x)
\end{bmatrix}\right\|_2,
\]
and
\[
\upsilon(x)=\frac52x^2,\qquad
\operatorname{dist}_2^2(-\nabla\upsilon(x),N_X(x))=25x^2=10\upsilon(x).
\]
At the boundary points \(x=\pm1\), the same identity holds because \(-\nabla\upsilon(x)\) points toward the interior of \(X\). Hence the feasibility--PL condition holds with \(\mu=5\) despite the failure of LICQ and MFCQ. Moreover, since \(X\) is bounded, the constraint values are also uniformly bounded, as required by the PG/GN reduction analysis.
\end{remark}

\section{Additional numerical implementation details}\label{app:numerical-details}
This appendix establishes the uniform parameter bound for the warm-started scalar regularization search described in \Cref{sec:Numerical_implementation}. 

\subsection{Warm-started scalar regularization}
\label{app:warm-start-search}
\begin{remark}[Boundedness of the scalar warm-started search]
\label{rem:warm-start-bound}
Under the assumptions giving $\tau_k\le\bar\tau$, this warm-started search also satisfies
\[
\underline M\le M_k\le\overline M
=2\max\{\underline M,M(\bar\tau)\},
\]
provided that the scalar main subproblem is solved exactly. Indeed, when $\operatorname{pred}_k>0$, \Cref{lem:basicproperties} gives
\[
\operatorname{pred}_k\ge M_k\|d^k\|_2^2,
\qquad
\zeta_k\ge1-\frac{L_{\nabla f}}{2M_k}.
\]
Consequently, $\zeta_k<0.1$ implies
$2M_k<10L_{\nabla f}/9\le\overline M$. The other carry-over branches do not increase $M_k$. Since $\widehat M_0=2\underline M\le\overline M$, induction bounds every initial trial $M_{k,0}$ by $\overline M$. If that initial trial is rejected, the trial immediately preceding acceptance is smaller than $M(\tau_k)$, because every value at least $M(\tau_k)$ is acceptable. Thus the accepted value is smaller than $2M(\tau_k)\le\overline M$. This proves both finite termination and the stated bound without requiring the Lipschitz constants in the implementation.
\end{remark}

\begin{footnotesize}
		\bibliographystyle{myplainurl}
		\bibliography{fsqp_verified}
	\end{footnotesize}

\end{document}